%% file: journal_nocomments.tex
\documentclass[a4paper,11pt]{article}

\usepackage{graphicx,reducecode}
\usepackage[bookmarks=false]{hyperref}
\usepackage{doi}
\usepackage[T1]{fontenc}
\usepackage[usenames,dvipsnames]{color}
\usepackage[dvipsnames]{xcolor}

\usepackage{amsthm,euler,amsmath,amssymb}
\usepackage{microtype}

\usepackage[authoryear]{natbib}

\title{
Boundary conditions for macroscale waves in an elastic system with microscale heterogeneity}

\author{Chen Chen\thanks{School of Mathematical Sciences, University of Adelaide, South Australia~5005, \textsc{Australia}.
\protect\url{mailto:chen.chen@adelaide.edu.au}}
\and
A. J. Roberts\thanks{School of Mathematical Sciences, University of Adelaide, South Australia~5005, \textsc{Australia}.
\protect\url{mailto:anthony.roberts@adelaide.edu.au}}
\and
J. E. Bunder\thanks{School of Mathematical Sciences, University of Adelaide, South Australia~5005, \textsc{Australia}.
\protect\url{mailto:Judith.Bunder@adelaide.edu.au}}
}

\date{\today}

\usepackage{pgfplots}
\usepackage{tikz}
\usepackage{pgf}

\usetikzlibrary{arrows,external}

\newcommand{\dd}[2]{\frac{\partial #1}{\partial #2}}
\newcommand{\DD}[2]{\frac{\partial^2 #1}{\partial #2^2}}

\newcommand{\abar}{\bar{a} }
\newcommand{\atilde}{\tilde{a} }
\newcommand{\Ord}{\mathcal{O}}

\newcommand{\lhs}{{\textsc{lhs}}}
\newcommand{\rhs}{{\textsc{rhs}}}
\newcommand{\ode}{{\textsc{ode}}}
\newcommand{\pde}{{\textsc{pde}}}

\usepackage{listings}
\usepackage{inconsolata}
\begin{document}

\maketitle

\begin{abstract}
Multiscale modelling aims to systematically construct macroscale models of materials with fine microscale structure.
However, macroscale boundary conditions are typically not systematically derived, but rely on heuristic arguments, potentially resulting in a macroscale model which fails to adequately capture the behaviour of the microscale system.
We derive the macroscale boundary conditions of the macroscale model for longitudinal wave propagation on a lattice with periodically varying density and elasticity.
We model the macroscale dynamics of the microscale Dirichlet, Robin-like, Cauchy-like and mixed boundary value problem.
Numerical experiments test the new methodology.
Our method of deriving boundary conditions significantly improves the accuracy of the macroscale models.
The methodology developed here can be adapted to a wide range of multiscale wave propagation problems.

\end{abstract}

\tableofcontents

\section{Introduction\label{sec:intro}}
This article develops a new rigorous method for deriving accurate boundary conditions for macroscale `homogenised' models of physical systems with significant microscale periodic heterogeneity.
Typically, for multiscale modelling of material mechanics with microscale structure, heuristic arguments are used to derive the macroscale boundary conditions.
For example, \citet[p.29]{MR2777986} ``expect'' a Dirichlet boundary condition for the macroscale model of a Dirichlet type boundary value problem in a composite material.
This expectation is reasonable when the scale separation between microscale and macroscale is very large.
This article significantly improves the accuracy of macroscale models, especially when the scale separation is not large, by deriving the macroscale boundary conditions systematically.
Examples in Sections~\ref{sub:First-iteration-1ch7} and~\ref{sec:otherbc} illustrate that the improvement from the usual heuristic boundary condition is by a factor comparable to the cell-to-domain ratio.

Many models of physical problems have counterintuitive initial conditions and boundary conditions.
For example, \citet{mercer1990acentre, Mercer:1994qf} built a rigorous macroscale model for shear dispersion in pipes in terms of cross-sectional averaged concentration.
They showed the effective initial and boundary conditions for the cross-sectional averaged concentration are not the same as the cross-sectional average of the initial and boundary conditions of the microscale system.
A long-known example is sound propagation in a closed-open tube, where, in the sound propagation model, the effective length of the tube is longer than the physical length \citep{Titze2001519, Parker1978245}.
The approach developed here provides a systematic and widely applicable rationale for quantifying such effects in the boundary conditions of macroscale models.

In multiscale modelling, \citet{Chen2014} proposed a cell-mapping methodology for deriving a macroscale model and corresponding boundary conditions for a one dimensional diffusion problem.
However, real applications of multiscale modelling techniques need the boundary condition methodology to be applicable to systems with multiple dimensions.
For example, a thin layer of metamaterial is often modelled as a two-dimensional periodic structure \citep[e.g.,][]{zbMATH05773414, Baron:13, 6167595, 6167594}.
The first step towards such multi-dimensional systems is to consider one-dimensional systems with an arbitrary number of strands.
Such strands may represent finite element\slash volume discretisations of the cross-section of a material in any number of dimensions.
This article develops previous methods \cite[e.g.]{Roberts92c} to the new domain of microscale heterogeneous spring-mass systems with arbitrary number of strands, arbitrary microscale longitudinal periodicities, and where the scale separation need not be very large.

\section{A spring-mass system with microscale structure\label{sec:micro}}
\begin{figure}
\centering\includegraphics[width=0.99\textwidth]{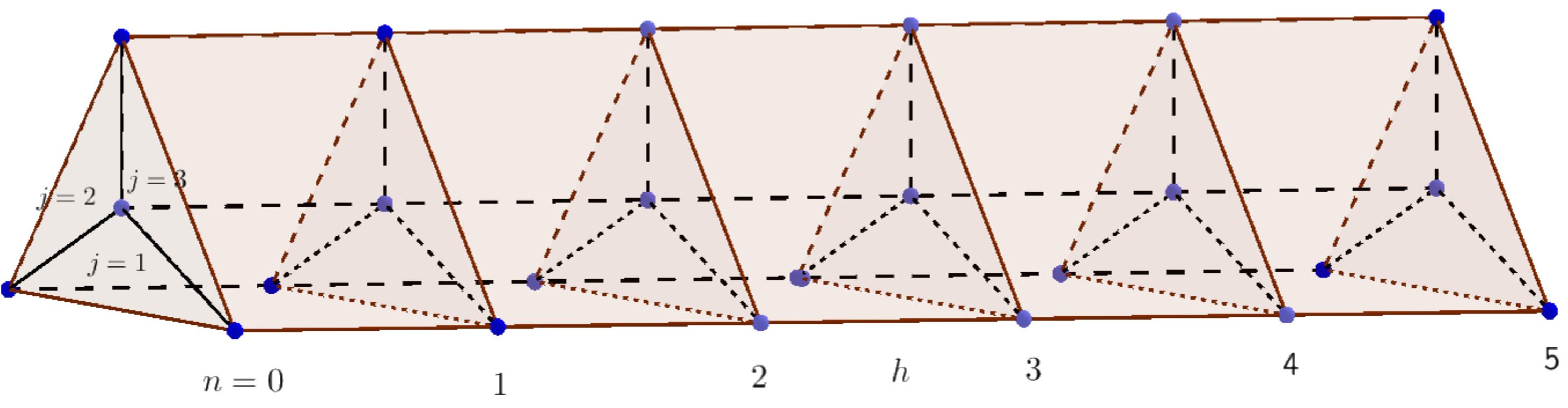}
\caption{A schematic diagram of a cross-sectional discretisation of a three-dimensional, relatively long and thin, object.
The illustrated prismatic, finite volumes notionally discretise the spatial dynamics into \(s=3\) `strands' and \(N=5\) longitudinal lattice intervals. The three prisms represent three different strands of the material.
\label{fig:schematic}}
\end{figure}

We model the  one-dimensional wave propagation along a metamaterial with `multiple strands'.
This article considers an $s$-strand spring-mass system on a microscale lattice material with longitudinal $p$-periodic density and elasticity, where $s$ and~$p$ are positive integers.
One way this model could arise is as a spatial discretisation (finite  element, finite volume, or multi-zone) of a three-dimensional heterogeneous elastic material as shown schematically in Figure~\ref{fig:schematic}.
For simplicity, we only consider linear interactions between the neighbouring elements.

The physically inspired variables are as follows.
Non-dimensionally, let \(x\)~measure position along the material, and let \(t\)~denote time.
We resolve the `cross-section' into \(s\)~strands indexed by $j=0,1,2,\ldots,s-1$ (Figure~\ref{fig:schematic}), and discretise the longitudinal dimension by a lattice with spacing~\(h\) and \(N\)~intervals.
Let the field~$u_{n,j}(t)$ be the displacement of a point-mass with nondimensional mass~\(h^3\rho_{n,j}\)%
, where $\rho_{n,j}$ is the effective density---on the \(j\)th~strand at longitudinal lattice point $n=0,1,\ldots,N$.
For this development of the methodology, we only model microscale vibrations in the longitudinal direction, and leave bending and torsion modes for further research.

The elasticities vary throughout the heterogeneous material.
On the $j$th~strand, the longitudinal elasticity between lattice point~$(n,j)$ and lattice point~$(n+1,j)$ is~$h\kappa_{n,j}$ (over a notional finite element cross-section of area~\(\Ord(h^2)\)).
The parameter~\(h\) denotes the equilibrium distance between any two adjacent masses in all direction (equivalently, the order of size of a spatial discretisation).
The cross elasticity~$h\kappa_{n,i,j}$  links two lattice points with longitudinal index~$n$ on the $i$th~strand
and $j$th~strand.
The elasticities~$\kappa_{n,i,j}$~and~$\kappa_{n,j}$ define the microstructure and are taken to vary with periodicity~$p$ in the longitudinal index~$n$.
Therefore, there are a total of $sp$~distinct longitudinal elasticities and $ps(s-1)/2$~distinct cross elasticities.
We assume no self elasticity, that is $\kappa_{n,j,j}=0$ for all~$n$ and~$j$, and that all the other elasticities and densities are strictly positive.

\begin{figure}
\centering
\includegraphics{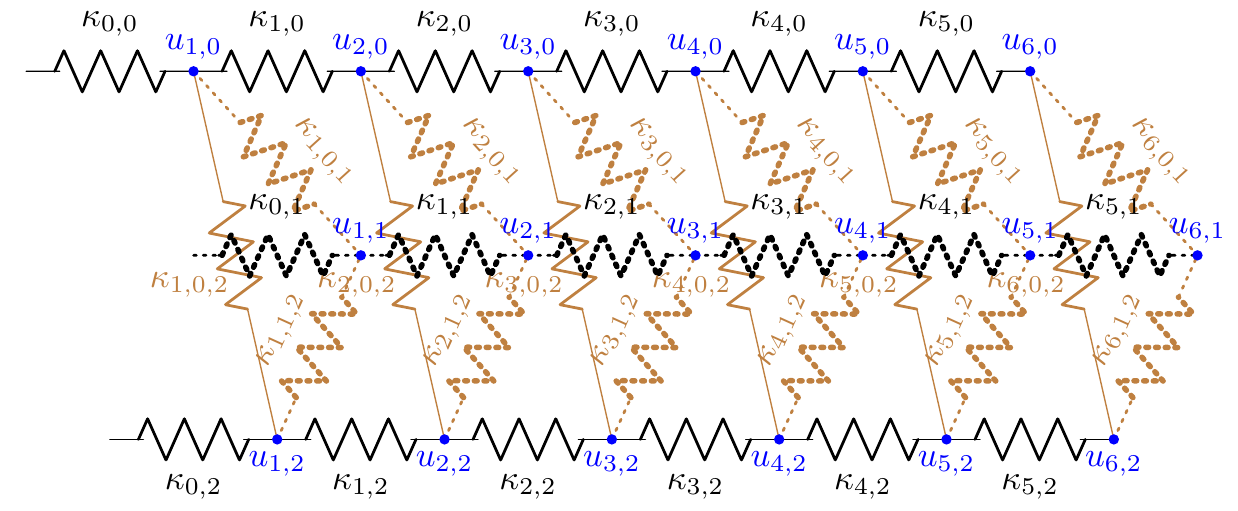}
\caption{A schematic diagram of a cell in a three-strand six-periodic spring mass system.
Parameter~$\kappa_{n,j}$ denotes the longitudinal elasticity between the $n$th and $(n+1)$th~masses on the $j$th~strand, whereas parameter~$\kappa_{n,i,j}$ denotes the cross elasticity between the $i$th and~$j$th strands.
The $n$th~mass on the $j$th~strand has mass~$h^3\rho_{n,j}$ and displacement~$u_{n,j}(t)$ from its equilibrium position.
\label{fig:s-layer-heterogeneous}}
\end{figure}

Define a cell to be the repetitive unit containing one period longitudinally and the entire cross-section: that is, a cell is formed by $s$~strands and $p$~consecutive masses on each strand. 
Hence the length of a cell is $ph$.
For example, Figure~\ref{fig:s-layer-heterogeneous} illustrates a cell of a six-periodic spring-mass system with three strands.
We model the longitudinal wave propagation in this multi-connected material.
An analogous material is studied for acoustic cloaking \citep{acousticsim, 1367-2630-9-3-045, 0022-3727-43-11-113001}.

By Hooke's law, the microscale dynamics on the $j$th~strand with $n$th~longitudinal index is governed by the coupled system of linear \ode{}s
\begin{eqnarray}
h^3\rho_{n,j}\DD{u_{n,j}}{t} & = & h\kappa_{n-1,j}\left[u_{n-1,j}(t)-u_{n,j}(t)\right]+h\kappa_{n,j}\left[u_{n+1,j}(t)-u_{n,j}(t)\right]\nonumber \\
 &  & {}+h\sum_{i=0}^{s-1}\kappa_{n,i,j}\left[u_{n,i}(t)-u_{n,j}(t)\right],\label{eq:sm}
\end{eqnarray}
for $n=1,2,\ldots,N-1$ and $j=0,1,\ldots,s-1$\,.
In the system of \ode{}s~\eqref{eq:sm}, the first two terms on the \rhs\ are due to the longitudinal springs and the last sum represents the inter-strand springs.
Sections~\ref{sec:micro}--\ref{sub:First-iteration-1ch7} assume that this spring-mass system has Dirichlet boundary conditions~$u_{0,j}(t)=b_{0,j}(t)$ and $u_{N,j}(t)=b_{N,j}(t)$ at the two ends of the longitudinal domain where \(b_{0,j}(t)\) and~\(b_{N,j}(t)\) are prescribed boundary values for strand $j=0$ to strand $j=s-1$\,.
Section~\ref{sec:otherbc} considers other microscale boundary conditions at the ends of the longitudinal domain.

\citet{Chen2014} provided explicit formulas for the macroscale model and boundary conditions for a two-strand two-periodic diffusion system.
Such explicit formulas are not feasible here because of the complexity of the algebra.
Instead, we develop, test and implement computer algebra (provided in Appendices~\ref{app:slowmanifoldmodel} and \ref{app:slowmanifoldbc}) which derives the macroscale model and macroscale boundary conditions for any given configuration of parameters.

Section~\ref{sec:ch7int} derives that the microscale heterogeneous system~\eqref{eq:sm} is modelled by the `homogenised' macroscale  wave \pde
\begin{equation}
\DD{U}{t}=\frac{\bar{\kappa}}{\bar{\rho}}\DD{U}{x}\,,\label{eq:chp7m}
\end{equation}
where constants~$\bar\kappa$ and $\bar\rho$ are the effective elasticities and densities, respectively. This sort of `homogenised' \pde\ model is well-known.
What is new is that we develop a methodology (Section~\ref{sub:First-iteration-1ch7}) to innovatively derive the corresponding macroscale Robin boundary conditions for homogenised macroscale \pde{}s such as~\eqref{eq:chp7m}:
\begin{equation}
U+d_{0}\dd{U}{x}  = B_{0}(t)\quad\text{at }x=0\,,
\quad\text{and}\quad
U+d_{L}\dd{U}{x} =  B_{L}(t)\quad\text{at }x=L\,.\label{eq:chp7bcf}
\end{equation}
The effective elasticity~$\bar{\kappa}$, the effective density~$\bar{\rho}$ and boundary condition coefficients~$d_{0}$ and~$d_{L}$ are complicated functions of the microscale elasticity and density distribution.
Previously, \citet{MR2777986}, \citet{PS08} and \citet{Bagdatli2015} presumed Dirichlet boundary conditions where $d_0,d_L=0$.
The numerical examples of Section~\ref{sec:ch7num} show that specific nonzero~$d_0$ and~$d_L$ for correctly chosen~\(B_0\) and~\(B_L\)  provide a significantly better macroscale model.

\section{Macroscale modelling in the interior \label{sec:ch7int}}
This section constructs the macroscale, homogenised, slow manifold \pde~\eqref{eq:chp7m}, valid in the domain interior, for the microscale spring-mass system~\eqref{eq:sm}.
Although the homogenised \pde~\eqref{eq:chp7m} is well-known, the approach here is not well-known, has not previously been applied to this class of problems, gives more powerful results, and better establishes the basis for systematically deriving improved boundary conditions.

\subsection{Fourier transform establishes the basis}

Recall that the material considered here is composed of cells repeated periodically where each cell has physical length~\(ph\) longitudinally.
Let $\nu=0,1,\ldots,\left\lfloor{N}/{p}\right\rfloor$ index cells  from the left boundary of the domain~$n=0$ to the right boundary \(n=N\)\,.
For example, $\nu=0$ is the left most cell, $\nu=1$ is the second left-most cell, and so on.
Also, define $m=n \bmod p\in\{0,1,\ldots,p-1\}$ which denotes the sub-cell longitudinal index of a mass.
Consequently, $u_{m+p\nu,j}(t)$ denotes the displacement of the mass with sub-cell longitudinal index~$m$ and strand index~$j$ within the $\nu$th~cell.

The modelling is best in Fourier space over cells, but best in physical space for the sub-cell structures.
Define~$\tilde{u}_{m,j}(k,t)$ to be a discrete Fourier transform over~\(\nu\) of displacement~$u_{m+p\nu,j}(t)$ over all cells across the whole domain.
Let~$S_k$ be a set of wavenumbers~$k$ such that eigenfunctions~$e^{ik(m+p\nu)h}$ are linearly independent and  for each~$m$ form a complete basis of~${u}_{m+p\nu,j}(t)$.
The Fourier expansion of the displacements is then
\begin{equation}
{u}_{m+p\nu,j}(t)=\sum_{k\in S_k}^{}
{\tilde{u}_{m,j}(k,t)}e^{ikh(m+p\nu)},
\label{eq:chap7defn}
\end{equation}
for sub-cell index \(m=0,\ldots,p-1\), strand index~$j=0,1,\ldots,s-1$\,, and  cells~\(\nu\).
Generalising the work of \citet{Chen2014}, substituting the Fourier expansions~\eqref{eq:chap7defn} into the \ode{}s~\eqref{eq:sm} and equating  coefficients of linearly independent functions~$e^{ik(m+p\nu)h}$ gives $sp$~coupled \ode{}s.
For each combination of~$m$ and~$j$, and all wavenumbers \(k\in S_k\),
\begin{eqnarray}
h^2\rho_{m,j}\DD{\tilde{u}_{m,j}}{t} & = & \kappa_{m-1,j}\left[\tilde{u}_{m-1,j}e^{-ikh}-\tilde{u}_{m,j}\right]+\kappa_{m,j}\left[\tilde{u}_{m+1,j}e^{ikh}-\tilde{u}_{m,j}\right]\nonumber \\
 &  & +\sum_{i=0}^{s-1}\kappa_{m,i,j}\left[\tilde u_{m,i}-\tilde u_{m,j}\right],\label{eq:smf}
\end{eqnarray}
where for convenience we define $\tilde{u}_{-1,j}(k,t):=\tilde{u}_{p-1,j}(k,t)$ and $\tilde{u}_{p,j}(k,t):=\tilde{u}_{0,j}(k,t)$.
The beauty of this Fourier transform is that in such linear problems each wavenumber~\(k\) is decoupled from all other wavenumbers, so we analyse them separately but inclusively with wavenumber~\(k\) as a parameter in the analysis.
Crucially, although the finite extent of the physical system limits relevant wavenumbers to \(k\in S_k\), we analyse the mathematical system for \emph{all} real wavenumbers.
Then the physically relevant subset of wavenumbers gives the physically relevant results.

To analyse the Fourier system~\eqref{eq:smf} we write it as a matrix-vector system for each real wavenumber~\(k\).
Define the vector of all Fourier-space fields across the spring-mass system
\[
\vec{u}(k,t):=\left(\tilde{u}_{0,0},\tilde{u}_{0,1},\tilde{u}_{0,2},\ldots,\tilde{u}_{0,s-1},\tilde{u}_{1,0},\tilde{u}_{1,1},\ldots,\tilde{u}_{1,s-1},\ldots,\tilde{u}_{p-1,s-1}\right),
\]
and define the corresponding $sp\times sp$ diagonal matrix of the densities
\[
B:=h^2\text{diag}\left(\rho_{0,0},\rho_{0,1},\ldots,\rho_{0,s-1},\rho_{1,0},\rho_{1,1},\ldots,\rho_{1,s-1},\ldots,\rho_{p-1,s-1}\right).
\]
Then we write the Fourier \ode{}s~\eqref{eq:smf} in the matrix-vector form
\begin{equation}
B\DD{\vec{u}}{t} = \mathcal{L}_{k}\vec{u}\,,\label{eq:smfd}
\end{equation}
where symmetric $sp\times sp$ matrix~$\mathcal{L}_{k}$ is almost a block cyclic
tridiagonal matrix (except for the nonzero top-right corner and bottom-left corner entries).
For longitudinal periodicity~$p=2$\,,
\begin{equation*}
\mathcal{L}_{k}  =  \begin{bmatrix}\mathcal{K}_{0} & K_{0}e^{ihk}+K_{1}e^{-ihk}\\
K_{0}e^{-ihk}+K_{1}e^{ihk} & \mathcal{K}_{1}
\end{bmatrix},
\end{equation*}
whereas for longitudinal periodicity~$p>2$\,,
\begin{equation*}
\mathcal{L}_k=\begin{bmatrix}
\mathcal{K}_0 & K_0e^{ihk} & 0&\cdots  &0& K_{p-1}e^{-ihk}\\
K_{0}e^{-ihk} & \mathcal{K}_1 & K_0e^{ihk} &0 &\cdots&0\\
0 & \ddots & \ddots & \ddots & \ddots & \vdots\\
\vdots & \ddots &\ddots & \ddots & \ddots &0\\
0&\cdots& 0 & K_{p-3}e^{-ihk} & \mathcal{K}_{p-2} &   K_{p-2}e^{ihk}\\
K_{p-1}e^{ihk} & 0& \cdots& 0 &K_{p-2}e^{-ihk} &\mathcal{K}_{p-1}
\end{bmatrix}.
\end{equation*}
Within the above matrix~$\mathcal{L}_{k}$, the off-diagonal blocks are the $s\times s$ symmetric matrices~$K_{m}=\text{diag}\left(\kappa_{m,0},\kappa_{m,1},\ldots,\kappa_{m,s-1}\right)$ which couples masses to their neighbours along each strand;
whereas the diagonal blocks are the $s\times s$ symmetric matrices  coupling the strands over a cross-section,
\[
\mathcal{K}_m=\begin{bmatrix}
-\kappa_{m,0}^+ & \kappa_{m,0,1} &\cdots& \cdots&\kappa_{m,0,s-1}\\
\kappa_{m,0,1}  & -\kappa_{m,1}^+ & \kappa_{m,1,2} &&\vdots \\
\vdots &\ddots &\ddots& \ddots&\vdots\\
\vdots&&\kappa_{m,s-3,s-2}  & -\kappa_{m,s-2}^+ & \kappa_{m,s-2,s-1}\\
\kappa_{m,0,s-1}&\cdots & \cdots & \kappa_{m,s-2,s-1} & -\kappa_{m,s-1}^+
\end{bmatrix}
\]
where the diagonals entries $\kappa^+_{m,i}=\kappa_{m-1,i}+\kappa_{m,i}+\sum_{j=0}^{s-1}\kappa_{m,i,j}$\,.
The factors of~\(e^{\pm ihk}\) appearing in~\(\mathcal L_k\) represent the Fourier transform coding of the coupling interactions between cells.

\subsection{There exists a slow manifold model}
\label{sec:tesmm}

As eigen modes corresponding to small wavenumber are most relevant to slow dynamics \citep{MEDiCS}, we base analysis on small wavenumber~\(k\).
Substitution verifies that a subspace of equilibria for the system~\eqref{eq:smfd}
occurs when the wavenumber $k=0$ and the Fourier displacements~\(\tilde u_{m,j}\) are equal at all lattice points.

To establish existence of a slow model we need to find the spectrum of the dynamics about each equilibria.
Because the system is linear, analysis of the equilibrium of zero displacement holds for all equilibria in the subspace of equilibria \citep{carr1981applications}.
Let $\mathcal{L}_{0}:=\left.\mathcal{L}_{k}\right|_{k=0}$ and
consider the base problem $B\partial{\vec{u}}/\partial{t}=\mathcal{L}_{0}\vec{u}$\,.
Seek a solution in the time harmonic form~$\vec{u}=e^{i\omega t}\vec{w}(t)$  where $\vec w(t)$~is an eigenvector corresponding to frequency~$\omega$.
Defining eigenvalues $\lambda:=\omega^2$ we hence need to solve the generalised eigen-problem
\begin{equation}
\lambda B\vec{{w}}+\mathcal L_{0}\vec{{w}}=0.
\label{eq:real1}
\end{equation}

\subsubsection{There exists eigenvalues of zero with multiplicity one}
Define vector~$\vec 1$ as a $sp$-dimensional constant vector with all components one.
Because the rows of the symmetric matrix~$\mathcal{L}_{0}$ sum to zero, direct substitution of the vector~$\vec w=\vec 1$ into the generalised eigen-problem~\eqref{eq:real1} shows that zero is an eigenvalue.
The eigenvalue of zero must have multiplicity of one because physically the spring mass system~\eqref{eq:sm} is only stationary when adjacent masses are spaced at the equilibrium distance~$h$.

\subsubsection{All nonzero eigenvalues are real and positive}

Both matrix~$\mathcal{L}_{0}$ and matrix~$B$ are real and symmetric.
Hence the Min-max Theorem \citep{Parlett1991373, teschlmathematical} give the real eigenvalues~$\lambda$ as the generalised Rayleigh quotient
\begin{equation}
\lambda=\frac{\vec{w}^{T}\left(-\mathcal{L}_{0}\right)\vec{w}}{\vec{w}^{T}B\vec{w}}\label{eq:chp7ray}
\end{equation}
which is a function of the corresponding eigenvector~$\vec w$.
Since all entries of the diagonal matrix~$B$ are positive,  $\vec w ^T B {\vec w}>0$\,.
It follows that eigenvalues~$\lambda\geq 0$ if and only if $\vec{w}^{T}\left(-\mathcal{L}_{0}\right)\vec{w} \geq 0$\,.
Recall that a matrix is diagonally dominant if, in each row, the magnitude of the diagonal entry is greater than or equal to the sum of the magnitude of the off diagonal entries \citep[p.881]{Kreyszig:2000:AEM:517902}.
By inspection, the magnitude of diagonals of matrix~$\left(-\mathcal{L}_{0}\right)$ is equal to the sum of the magnitudes of off diagonal entries in matrix~$\left(-\mathcal{L}_{0}\right)$.
Hence the symmetric matrix~$\left(-\mathcal{L}_{0}\right)$ is diagonally dominant.
Also, the diagonal entries of matrix~$\left(-\mathcal{L}_{0}\right)$ are all positive.
A symmetric diagonally dominant matrix with nonnegative diagonals is positive semidefinite \citep{definitesemi}.
Hence the inner product~$\vec{w}^{T}\left(-\mathcal{L}_{0}\right)\vec{w}\geq0$\,.
Substituting $\vec{w}^{T}\left(-\mathcal{L}_{0}\right)\vec{w}\geq0$ and $\vec w ^T B {\vec w}>0$ into the generalised Rayleigh quotient~\eqref{eq:chp7ray}
implies the eigenvalues of the generalised eigen-problem~\eqref{eq:real1}
are all non-negative.

\subsubsection{All positive eigenvalues are separated from zero}

There must exist a spectral gap between the positive eigenvalues and the zero eigenvalue because the base system~\eqref{eq:real1} is finite dimensional and has only a finite number of eigenvalues.
Hence all positive eigenvalues are separated from the zero eigenvalue by a finite gap.

\subsubsection{The spectrum of the linearised spring-mass system}

Recall that from the base eigen-problem~\eqref{eq:real1}, the frequencies for the spring-mass system are~$\omega=\pm\sqrt{\lambda}$\,.
Hence the zero eigenvalue of the generalised eigen-problem~\eqref{eq:real1} corresponds to repeated zero frequencies of the linearised spring-mass system~$B\partial{\vec{u}}/\partial{t}=\mathcal{L}_{0}\vec{u}$\,.
All positive eigenvalues of the generalised eigen-problem~\eqref{eq:real1} corresponds to two purely oscillatory modes of the spring-mass system with frequency~$\omega=\pm\sqrt{\lambda}$ well separated from zero.
Hence sub-centre manifold theory by \citet[Thm.~7.1]{0577.34039} assures us of the existence of a slow manifold macroscale model to system~\eqref{eq:smfd} in a finite neighbourhood of the subspace of equilibria with wavenumber $k=0$\,.

\subsection{Derive a slow manifold model for low wavenumber}

The next step is to construct the slow manifold model as a power series in the wavenumber~$k$ (because Section~\ref{sec:tesmm} establishes existence for small wavenumber~\(k\)).
By Taylor series expansion, $e^{\pm ikh}=1\pm ikh-\frac{1}{2}h^2k^{2}+\mathcal{O}(k^{3})$.
Thus we approximate system~\eqref{eq:smfd} for small wavenumber~$k$ with
\begin{equation}
B\DD{\vec{u}}{t}=\left.\mathcal{L}_{k}\right|_{e^{\pm ikh}=1\pm ikh-\frac{1}{2}h^2k^{2}}\vec{u}+\mathcal{O}(k^{3}).\label{eq:taylapp2}
\end{equation}
Defining the velocity Fourier variables $\vec{v}(k,t):=\dd{\vec{u}}{t}$, we transform system~\eqref{eq:taylapp2}
into dynamical system form
\begin{equation}
\dd{}{t}\begin{bmatrix}\vec{u}\\
B\vec{v}
\end{bmatrix}=\begin{bmatrix}0 & I\\
\left.\mathcal{L}_{k}\right|_{e^{\pm ik}=1\pm ikh-\frac{1}{2}h^2k^{2}} & 0
\end{bmatrix}\begin{bmatrix}\vec{u}\\
\vec{v}
\end{bmatrix}+\mathcal{O}(k^{3}).\label{eq:taylap}
\end{equation}
We aim to derive a slow manifold model to errors~$\Ord(k^{3})$ for the dynamical system~\eqref{eq:taylap}.
Due to the complex structure of matrix~$\mathcal{L}_{k}$,
it is algebraically extremely complicated to manually derive a slow manifold model for the system~\eqref{eq:taylap}.
Instead, we invoke computer algebra provided in Appendices~\ref{app:slowmanifoldmodel} and \ref{app:slowmanifoldbc} to construct a slow manifold model.

\subsection{Amplitude definition and initial approximation \label{sub:ini}}

Define the macroscale Fourier space amplitude to be the cell-average
\begin{equation}
\tilde{U}(k,t):=\frac{1}{sp}\sum_{j=0}^{s-1}\sum_{m=0}^{p-1}\tilde{u}_{m,j}(k,t).\label{eq:amplitudeu}
\end{equation}
We need an extra macroscale variable to parameterise the slow manifold. The zero frequency given by the generalised eigenvalue problem~\eqref{eq:real1} is of multiplicity two. Thus, the slow manifold is two dimensional for each wavenumber~\(k\). This implies we need another macroscale variable.
Consequently, define the macroscale Fourier velocity as
\begin{equation}
\tilde{V}(k,t):=\frac{1}{sp}\sum_{j=0}^{s-1}\sum_{m=0}^{p-1}\dd{}{t}\tilde{u}_{m,j}(k,t).
\label{eq:amplitudev}
\end{equation}
Then an initial approximation to the slow manifold for the perturbed linear system~$B{\partial\vec{u}}/{\partial t}=\mathcal{L}_{0}\vec{u}$, in the Fourier slow subspace, is
\begin{equation}
\vec{u}=\tilde{U}(k,t)\vec{1}\quad\text{and}\quad\vec{v}=\tilde{V}(k,t)\vec{1}\,,\label{eq:inisch7}
\end{equation}
with evolution
\begin{equation}
\dd{\tilde{U}}{t}=\tilde V\quad\text{and}\quad\dd{\tilde{V}}{t}=\tilde g(\tilde U,k)= 0\,.\label{eq:iniech7}
\end{equation}

\subsection{Macroscale model and examples \label{sub:slowm}}

We implement computer algebra provided in Appendix~\ref{app:slowmanifoldmodel} to construct the macroscale model of microscale wave \ode{}s~\eqref{eq:smfd}.

The algorithm of Appendix~\ref{app:slowmanifoldmodel}  firstly (lines~9--50)  encodes the wave system and its initial approximations as given in Section~\ref{sub:ini}.
Secondly, the algorithm (lines~51--68) iteratively seeks better approximations for the slow manifold~$\vec{u}_\text{new}=\vec{u}_\text{old}+\hat{\vec u}(\tilde U, \tilde V)$ and~$\vec{v}_\text{new}=\vec{v}_\text{old}+\hat{\vec v}(\tilde U, \tilde V)$, where vectors~$\hat{\vec{u}}(\tilde{U},\tilde{V})$ and~$\hat{\vec{v}}(\tilde{U},\tilde{V})$ are the corrections to the shape of the slow manifold~\eqref{eq:inisch7}.
The algorithm computes the residuals of the governing \pde{}s~\eqref{eq:taylapp2} and uses the residuals to derive these corrections.
The iteration also uses the residuals to correct the right-hand side~\(g\) of the evolution~\eqref{eq:iniech7} on the slow manifold.
Let $\vec 0$ denote the $sp$-dimensional vector with all components zero.
The two amplitude definitions~\eqref{eq:amplitudeu} and~\eqref{eq:amplitudev} give two constraints for the corrections
of the slow manifold
\begin{equation}
\vec 1^T\hat{\vec{u}}(\tilde{U},\tilde{V})=\vec 0\quad\text{and}\quad
\vec 1^T\hat{\vec{v}}(\tilde{U},\tilde{V})=\vec 0\,;\label{eq:ch7amp}
\end{equation}
line~52 encodes these two constraints.
The algorithm terminates the iteration when the computed residual is smaller than the truncation order specified in line~21.
Finally, the algorithm outputs the slow manifold and its evolution into a text file.

For example, the algorithm of Appendix~\ref{app:slowmanifoldmodel} derives that a simple two-strand two-periodic spring-mass system~\eqref{eq:sm} has macroscale model~\eqref{eq:chp7m} with effective density as the average
\begin{equation*}
\bar{\rho}=\tfrac{1}{4}\sum_{i=0}^{1}\sum_{j=0}^{1}\rho_{i,j},
\end{equation*}
and the effective elasticity
\begin{equation*}
\bar{\kappa}=\frac{\kappa_{0}\kappa_{1}\left(\kappa_{0,1}+\kappa_{0,0}\right)\left(\kappa_{1,1}+\kappa_{1,0}\right)+\left(\kappa_{1,0,1}+\kappa_{0,0,1}\right)\kappa_{0,0}\kappa_{0,1}\kappa_{1,0}\kappa_{1,1}\sum_{j=0}^{1}\sum_{i=0}^{1}1/\kappa_{i,j}}{\kappa_{0}\kappa_{1}(\kappa_{0,0}+\kappa_{0,1}+\kappa_{1,0}+\kappa_{1,1})+(\kappa_{1,0,1}+\kappa_{0,0,1})(\kappa_{1,1}+\kappa_{0,1})(\kappa_{1,0}+\kappa_{0,0})}.
\end{equation*}
Although the iterative code of Appendix~\ref{app:slowmanifoldmodel} can derive analytic expressions such as these, symbolic results for $s$-strand $p$-periodic systems are usually not helpful because the expressions are extremely complicated.
In these cases, for this article, we substitute sample numerical values of the elasticities and densities into analytical formulae.

For systems that has a higher number of strand and periodicity, analytical derivations are not feasible and we implement numerical computations. 
For example, for a five-strand and ten-periodic system, set the prescribed elasticity distributions to
\begin{equation}\label{eq:ela}
\kappa_{n,j}=\frac{1}{1+A_j\cos(4.6nh+\phi_j)} \quad \text{and} \quad \kappa_{n,j_1,j_2}=\frac{1}{1+A_{j_1,j_2}'\cos(4.6nh+\phi_{j_1,j_2}')},
\end{equation}
and the prescribed density distributions are
\begin{equation}\label{eq:den}
\rho_{n,j}={1+B_j\sin(4.6nh+\varphi_j)},
\end{equation}
where table~\ref{tab:coe} gives the coefficients~$A_j$, $A_{j_1,j_2}'$, $B_j$, $\phi_j$, $\varphi_j$ and~$\phi_{j_1,j_2}'$.
\begin{table}
\caption{Coefficients for the elasticity and density distributions in the five-strand
ten-periodic example.\label{tab:coe}}
\begin{equation*}
\begin{array}{l|rrrrr}
j & 1 & 2 & 3 & 4 & 5 \\
\hline
A_{j} & 0.929  & 0.776  & 0.487  & 0.436  & 0.447 \\
B_{j} & 0.963 & 0.547 & 0.521 & 0.231 & 0.489 \\
\phi_{j} & -1.217  & 0.053  & 0.068  & 1.996  & 1.852 \\
\varphi_{j} & 0.779  & 1.126  & -0.656  & -0.833  & 3.066 \\
A_{1,j} & 0  & 0.939  & 0.208  & 0.195  & 0.311 \\
A_{2,j} & 0.939  & 0  & 0.301  & 0.226  & 0.923 \\
A_{3,j} & 0.208  & 0.301  & 0  & 0.171  & 0.430 \\
A_{4,j} & 0.195  & 0.226  & 0.171  & 0  & 0.185 \\
A_{5,j} & 0.311  & 0.923  & 0.430  & 0.185  & 0 \\
\phi_{1,j} & 0  & 0.596  & -2.404  & -2.604  & 1.447 \\
\phi_{2,j} & 0.596  & 0  & -1.278  & -1.492  & -0.0720 \\
\phi_{3,j} & -2.404  & -1.278  & 0  & 1.891  & 0.493 \\
\phi_{4,j} & -2.604  & -1.492  & 1.891  & 0  & -1.651 \\
\phi_{5,j} & 1.447  & -0.0720  & 0.493  & -1.651  & 0 \\
\end{array}
\end{equation*}
\end{table}
The computer algebra of Appendix~\ref{app:slowmanifoldmodel} derives the slow manifold of this example as
\begin{equation}\label{eq:ch7ex}
\vec{u}(\tilde U,t)=(\vec 1+ik\vec \alpha+k^2\vec\beta)\tilde U(k,t)+\Ord(k^3),
\end{equation}
where  the constant vectors~$\vec\alpha$ and~$\vec\beta$ (listed in Appendix~\ref{app:slowmanifold}) have means of zeros and standard deviations of $0.46$ and $0.64$, respectively. 
The evolution on this slow manifold is
\begin{equation}
\DD{\tilde U}{t}=-1.176\,k^2\tilde U(k,t)+\Ord(k^3). \label{eq:ch7evoe}
\end{equation}
As discussed by \citet{MEDiCS}, this macroscale model is valid for small but finite wavenumber, $|k|\leq K$ for some~\(K\), hence it is a model for the slowly varying dynamics of the microscale heterogeneous elastic spring-mass system~\eqref{eq:smfd}.

Although ideally the spring-mass system is not dissipative, damping or energy radiation is always present in real life applications \citep{1997theoretical}.
Consequently, oscillatory modes with wavenumbers~$k>|K|$ decay, especially those of high wavenumber.
As a result, the slow manifold dynamics captured by this modelling is observed in real physical systems.

Define the macroscale spatial variable~$x$ to be the distance from the left boundary $n=0$\,.
For example, at the $n$th~lattice points,~$x=nh$\,.
Taking the inverse Fourier transform of the evolution~\eqref{eq:ch7evoe} gives the physical domain, homogenised, wave \pde
\begin{equation}\label{eq:chp7macromodel}
\DD{U(x,t)}{t}=1.176\DD{U(x,t)}{x}+\Ord(\partial^3_x).
\end{equation}
The construction of this section may be extended to derive higher-order \pde\ models, but this classic wave \pde\ is sufficient to show the improvements upon using good boundary conditions.
The next section establishes the methodology to derive the macroscale boundary condition of such macroscale \pde{}s.

\section{Derive boundary conditions for slow manifold
model\label{sub:First-iteration-1ch7}}

This section details the innovations needed to derive the macroscale Robin boundary conditions~\eqref{eq:chp7bcf} for the macroscale homogenised \pde~\eqref{eq:chp7m} of the microscale spring-mass system~\eqref{eq:sm}.

The derivation of the macroscale boundary conditions~\eqref{eq:chp7bcf} for the macroscale
\pde~\eqref{eq:chp7m} employs a local analysis in the longitudinal space.
This section derives the coefficients~$d_{0}$, $d_{L}$, $B_0$ and~$B_L$ in the macroscale Robin boundary conditions~\eqref{eq:chp7bcf} using the density and elasticity microscale details in boundary layers and towards the interior.
The derivation of the macroscale boundary condition
coefficients involves solving an $sp\times sp$
linear system, solving a $2s$~dimensional eigenvalue problem, and
finding the null space of a $\left(s+2\right)\times\left(s+1\right)$
matrix.

For large numbers of strands~$s$ it is infeasible to do these
steps manually, but they are computationally inexpensive to derive either by numerical computation or by computer algebra.
If the coefficients of the microscale system are numerically specified, then the computer algebra program of Appendix~\ref{app:slowmanifoldbc} derives numerical values for the coefficients of the macroscale boundary conditions~\eqref{eq:chp7bcf}.
Alternatively, if the coefficients are not numerically specified, then the computer algebra program of Appendix~\ref{app:slowmanifoldbc} symbolically derives the macroscale boundary conditions~\eqref{eq:chp7bcf} as algebraic expressions in terms of the density~$\rho_{n,j}$ and elasticity constants~$\kappa_{n,j}$ and~$\kappa_{n,j_{1},j_{2}}$.
This section gives explicit formulas for the general two-strand two-periodic system as an algebraically accessible example.

\subsection{Left-end macroscale boundary condition\label{sub:left}}

Since we seek a slowly varying macroscale model for the microscale dynamics, the spatial structures of the microscale system~\eqref{eq:sm} is responsible for boundary layers in the microscale solution.
We assume time derivatives are negligible when deriving macroscale boundary conditions for the slowly varying macroscale model.
This assumption simplifies our derivation of the macroscale boundary conditions, is plausible because spatial structures dominant the boundary layers, and the assumption is tested to be reasonable (Sections~\ref{sec:ch7num2} and~\ref{sec:ch7num}).
To proceed, note that any quasi-steady state of the microscale spring-mass system~\eqref{eq:sm} satisfies
\begin{eqnarray}
0 &=&  \kappa_{n-1,j}\left[u_{n-1,j}(t)-u_{n,j}(t)\right]+\kappa_{n,j}\left[u_{n+1,j}(t)-u_{n,j}(t)\right]\nonumber\\
&& {}+\sum_{i=0}^{s-1}\kappa_{n,i,j}\left[u_{n,i}(t)-u_{n,j}(t)\right],
\quad\text{for all }n,j,
\qquad\label{eq:smss}
\end{eqnarray}
with Dirichlet boundary conditions~$u_{0,j}(t)=b_{0,j}(t)$ and $u_{N,j}(t)=b_{N,j}(t)$.
These quasi-steady state equations form an algebraic linear system.

We write the linear equations~\eqref{eq:smss} into a matrix-vector form. Define the vectors
\begin{eqnarray*}&&
\vec{u}_{0}  :=  (u_{0,0},
\ldots,
u_{0,s-1},
u_{1,0},
\ldots,
u_{1,s-1}
)
\quad\text{and}\\&&
\vec{u}_{1}  :=  (u_{p,0},
\ldots,
u_{p,s-1},
u_{p+1,0},
\ldots,
u_{p+1,s-1}
).
\end{eqnarray*}
These two vectors contain the displacements of the leftmost $2s$~masses in each of the first two cells.
We put the first $sp$~rows
of the quasi-steady state spring-mass system~\eqref{eq:smss} into a matrix-vector form.
As an example, for a two-strand two-periodic system, the first $sp=4$~rows of~\eqref{eq:smss} are
\begin{equation}\label{eq:diff}
\vec 0=
\begin{bmatrix}\kappa_{0,0} & 0 & -\kappa_{1,0}^{+} & \kappa_{0,1,1} & \kappa_{1,0} & 0 & 0 & 0\\
0 & \kappa_{0,1} & \kappa_{0,1,1} & -\kappa_{1,1}^{+} & 0 & \kappa_{1,1} & 0 & 0\\
0 & 0 & \kappa_{1,0} & 0 & -\kappa_{0,0}^{+} & \kappa_{0,1,0} & \kappa_{0,0} & 0\\
0 & 0 & 0 & \kappa_{1,1} & \kappa_{0,1,0} & -\kappa_{0,1}^{+} & 0 & \kappa_{0,1}
\end{bmatrix}\begin{bmatrix}
\vec{u}_0\\
\vec{u}_1
\end{bmatrix},
\end{equation}
where
\begin{align*} &
\kappa_{1,0}^{+} :=\kappa_{0,0}+\kappa_{1,0}+\kappa_{0,1,1}, &&
\kappa_{1,1}^{+} := \kappa_{0,1}+\kappa_{1,1}+\kappa_{0,1,1},\nonumber \\ &
\kappa_{0,0}^{+} := \kappa_{1,0}+\kappa_{0,0}+\kappa_{0,1,0}, &&
\kappa_{0,1}^{+}:= \kappa_{1,1}+\kappa_{0,1}+\kappa_{0,1,0}.
\end{align*}
For systems whose periodicity is greater than two, the vector on the \rhs\ of equation~\eqref{eq:diff} is of the form
\begin{math}
\begin{pmatrix}\vec u_0,u_{2,0},
\ldots,
u_{2,s-1},
\ldots,
u_{p-1,s-1},\vec u_1 \end{pmatrix}.
\end{math}
We then arrange the equilibrium equations into a cell mapping form~$\vec u_1=T \vec u_0$ where the $2s \times 2s$ matrix~$T$ is derived from the first $sp$~equations in the quasi-steady state spring-mass system~\eqref{eq:smss}.
For example, the two-periodic two-strand system has the matrix
\begin{eqnarray}\label{eq:matrixTch6}
T&:=&-\begin{bmatrix}\kappa_{1,0} & 0 & 0 & 0\\
0 & \kappa_{1,1} & 0 & 0\\
-\kappa_{1,0}-\kappa_{0,0}-\kappa_{0,1,0} & \kappa_{0,1,0} & \kappa_{0,0} & 0\\
\kappa_{0,1,0} & -\kappa_{1,1}-\kappa_{0,1}-\kappa_{0,1,0} & 0 & \kappa_{0,1}
\end{bmatrix}^{-1}
\nonumber\\
&&{}\times\begin{bmatrix}\kappa_{0,0} & 0 & -\kappa_{0,0}-\kappa_{1,0}-\kappa_{0,1,1} & \kappa_{0,1,1}\\
0 & \kappa_{0,1} & \kappa_{0,1,1} & -\kappa_{0,1}-\kappa_{1,1}-\kappa_{0,1,1}\\
0 & 0 & \kappa_{1,0} & 0\\
0 & 0 & 0 & \kappa_{1,1}
\end{bmatrix}.\qquad
\end{eqnarray}
We compute the $2s$ eigenvalues and eigenvectors of the matrix~$T$.
Numerical results demonstrate there is always an eigenvalue of one (repeated twice)
with one eigenvector~$\vec{1}$ and one generalised eigenvector.\footnote{More than five hundred realisations numerically confirm this property.  The elasticity and density are randomly selected within zero and a hundred.  The periodicity is selected between two and fifteen and the number of strands is selected to be within one and ten. Further research aims to establish this property of the eigen-spectrum by analysis.}
The span of these two eigenvectors are the centre subspace of matrix~$T$.
There are also $s-1$ real eigenvalues greater than one which correspond to exponentially growing modes, and $s-1$ real eigenvalues strictly between zero and one which correspond to exponentially decaying modes.
Let~$\mu_i$ for $i=0,1,\ldots,2s-1$ denote the eigenvalues of the matrix~$T$.
For definiteness, order the eigenvalues so that $\mu_1\leq\mu_2\leq\cdots\leq\mu_{2s}$; consequently $\mu_{s}=\mu_{s+1}=1$\,.
For example, in the two-periodic two-strand problem, the mapping matrix~$T$ defined in~\eqref{eq:matrixTch6} has repeated eigenvalue~$\mu_2=\mu_3=1\,$, an eigenvalue~$1>\mu_1>0$ and an eigenvalue~$\mu_4>1$\,.
\citet{opac-b1078099} discussed in detail the trichotomy of the eigenspace for maps and discrete dynamical systems.  This article views the longitudinal spatial structure of the material as a `time'-like evolution as we move away from the left boundary into the interior, and derives the macroscale boundary conditions as a consequence.

We write the displacement vector~$\vec u_0$ as a linear combination of the eigenvectors of the mapping matrix~$T$.
On the boundary~$x=0$, we are only interested in the centre modes and decaying modes.
This restriction is because if any growing mode is present in the left boundary layer, then it would be exponentially large at the right boundary.
This largeness is not physically acceptable and hence there must be no exponentially growing modes in the left boundary layer.
As there is no exponentially growing unstable modes, the boundary layer solution near \(x=0\) must be a linear combination of the~$s+1$ eigenvectors corresponding to the centre-stable modes.
Hence at the left-end of the first cell,
\begin{equation}
\vec{u}_{0}  =   \begin{pmatrix}b_{0,0},
\ldots,
b_{0,s-1},
b_{1,0},
\ldots,
b_{1,s-1}
\end{pmatrix}=\sum_{i=1}^{s+1}c_{i}\vec{v}_{i},\label{eq:lincom}
\end{equation}
where $\vec v_i$ where $i=1,2,\ldots,s,s+2,\ldots,2s$ are eigenvectors corresponding to eigenvalue~$\mu_i$,
and $\vec v_{s+1}$~is a generalised eigenvector corresponding to eigenvalue one.
The~$s$ microscale Dirichlet boundary conditions give~$s$ constraints on the boundary with~$s+1$ unknown coefficients~$c_i$.

Our new approach finds the macroscale boundary conditions by a projection onto the centre-stable manifold.
Recall from Section~\ref{sec:ch7int} that the macroscale displacement~$U(x,t)$ models the dynamics of the microscale system~\eqref{eq:sm} in the domain interior.
So we specify the macroscale boundary conditions which gives an accurate interior macroscale model.
The exponential trichotomy of manifolds assures us the solution in the left boundary layer is on the centre-stable manifold.
When we move away from the left boundary into the domain interior, the stable modes decay exponentially quickly and become negligible.
Hence the macroscale displacement~$U(x,t)$, which aims to model the domain interior, is a linear combination of the basis of the centre modes.
Project the~$u_{0,j}(t)$ onto a two dimensional centre subspace by setting the $s-1$ coefficients of stable eigenvectors to zero: this projection corresponds to the exponential decay of the stable modes within the boundary layer \cite[]{Roberts92c}.
Then we use the results of the projection and definitions~\eqref{eq:amplitudeu} and~\eqref{eq:amplitudev} to compute macroscale variables~$U(0,t)$ and~$V(0,t)$ on the centre subspace.
The macroscale amplitudes are initially defined in the domain interior.
An extrapolation defines~$U(x,t)$ and~$V(x,t)$ on the boundaries.
These definitions of~$U(0,t)$ and flux~$V(0,t)$ give another two constraints on the boundary values of the macroscale fields.
With the previous constraints from microscale boundary values~$u_{0,j}(t)$, there are in total $s+2$ equations and $s+1$ unknowns and they form a $(s+2)\times(s+1)$ linear system.
Consequently, the right-hand sides must satisfy a solvability condition in order for a solution to exist: this solvability condition is the effective macroscale boundary condition.

The solvability condition is that the right-hans side must be orthogonal to the null space of the transpose matrix.
The transpose, being $(s+1)\times(s+2)$, generically poses a one dimensional null space spanned by one $s+2$ dimensional basis vector.
The dot product of this basis vector and the vector
\begin{equation*}
\left(b_{0,0}(t),b_{0,1}(t),\ldots,b_{0,s-1}(t),\left.U\right|_{x=0},\left.\dd{U}{x}\right|_{x=0}\right),
\end{equation*}
which contains all the $s$~microscale boundary values~$u_{0,j}(t)$ and the amplitude definition of the macroscale boundary values~$U(0,t)$ and flux~$V(0,t)$, gives the macroscale boundary condition at $x=0$\,.
This procedure is implemented in general by the computer algebra program of Appendix~\ref{app:slowmanifoldbc}.

\begin{figure}
\centering \includegraphics{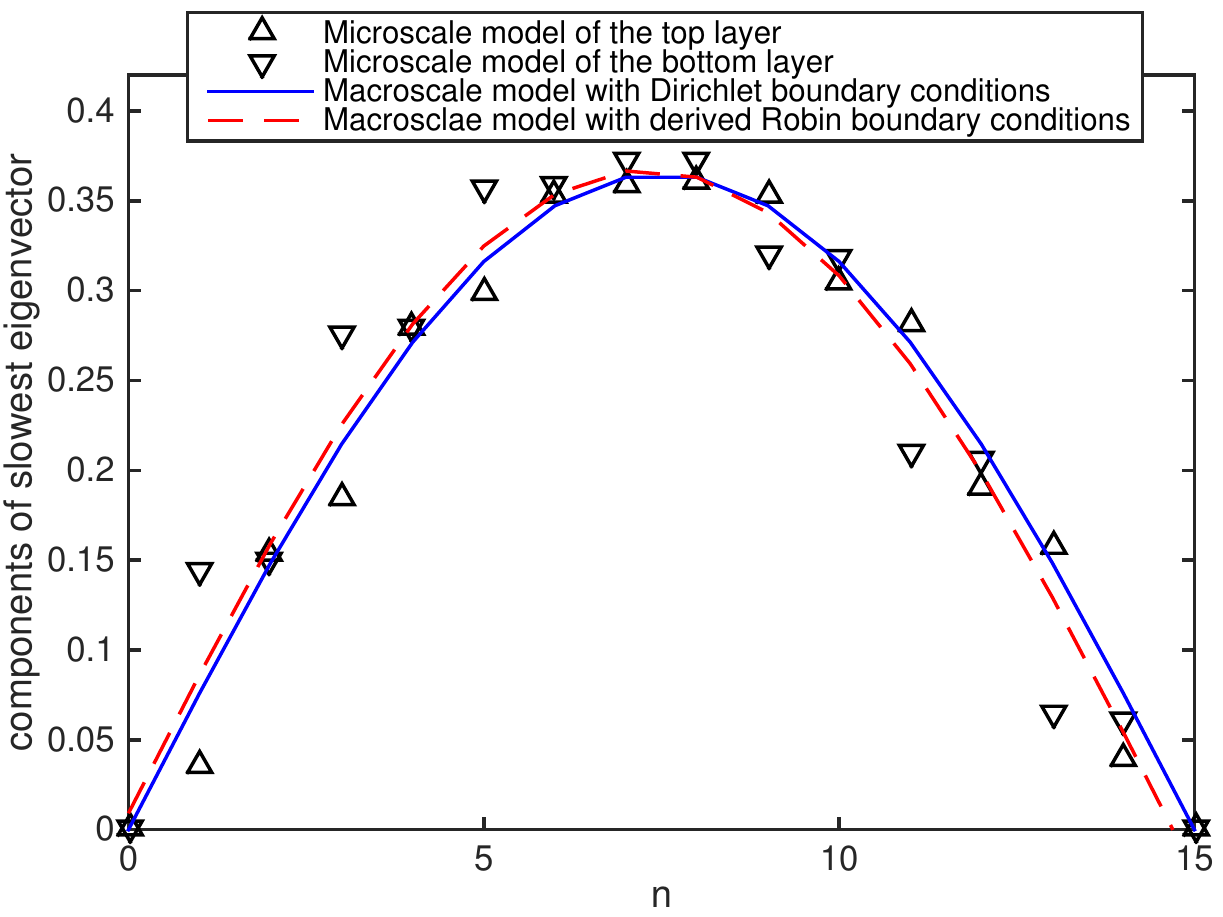}
\caption{The slowest eigenvector components of the spring-mass system~\eqref{eq:sm} with a two strands and periodicity two.
The longitudinal elasticities are $\kappa_{0,0}=2, \kappa_{0,1}=0.5, \kappa_{1,0}=0.1, \kappa_{1,1}=5$\,, the cross elasticities are $\kappa_{0}=1,\kappa_1=0.1$  and the densities are $\rho_{0,0}=1, \rho_{0,1}=2, \rho_{1,0}=4, \rho_{1,1}=0.5$ respectively.
The triangles are the slowest eigenvector of the microscale system.
The red dotted line is our macroscale model~\eqref{eq:chp7m} with derived Robin boundary conditions~\eqref{eq:paper1bc}
that correctly produce the solution in the interior of the domain.
The blue line is the macroscale model obtained with classic Dirichlet boundary conditions $U(0,t)=U(L,t)=0$ that are not quite as good.\label{fig:ch6num}}
\end{figure}

We illustrate the method with two strands and periodicity two as an example.
Let~$v_{ij}$ denote the $j$th~components of the eigenvector~$\vec v_i$.
In the case of two-periodic two-strand, the eigenvector $\vec v_{2}=(1,1,1,1)$ corresponds to the eigenvalue of one.
Physically, the eigenvector~$(1,1,1,1)$ describes the
equilibrium when no sub-cell elastic deformation is present.
The $4\times 3$ linear system describing boundary constraints is
\begin{equation}
\begin{bmatrix}v_{11} & 1 & v_{31}\\
v_{12} & 1 & v_{32}\\
0 & 1 & \tfrac{1}{4}\left(\vec{v}_{3}\cdot\vec{1}-1\right)\\
0 & 0 & \frac{1}{2h}
\end{bmatrix}\begin{bmatrix}c_{1}\\
c_{2}\\
c_{3}
\end{bmatrix}=\begin{bmatrix}b_{0,0}\\
b_{0,1}\\
\left.U\right|_{x=0}\\
\left.\dd{U}{x}\right|_{x=0}
\end{bmatrix}.\label{eq:twolayer-UdUdx}
\end{equation}
The first two rows of this matrix system come from the linear combination~\eqref{eq:lincom} and the last two rows come from the definitions~\eqref{eq:amplitudeu} and \eqref{eq:amplitudev}.
The vector spanning the null space of the transpose of the $4\times3$~matrix in this linear system is
\begin{equation*}
\vec w=\begin{bmatrix}
v_{12}/(-v_{12}+v_{11})\\
-v_{11}/(-v_{12}+v_{11})\\
1\\
-2h\left(\frac{v_{12}v_{31}}{-v_{12}+v_{11}}-\frac{v_{11}v_{32}}{-v_{12}+v_{11}}+\tfrac{1}{4}\big(\vec{v}_3\cdot\vec{1}-1\big)\right)
\end{bmatrix}.
\end{equation*}
The vector~$\vec w$ always exists as the matrix in equation~\eqref{eq:twolayer-UdUdx} is $4\times3$ so the null space is never empty.
Premultiply~\eqref{eq:twolayer-UdUdx} by~$\vec w^T$
to obtain zero on the \textsc{lhs}.
Thus the \textsc{rhs} of~\eqref{eq:twolayer-UdUdx}, premultiplied by~$\vec w^T$, provides the required boundary condition at $x=0$ of
\begin{equation}
U-2h\left[\frac{v_{12}v_{31}-v_{11}v_{32}}{-v_{12}+v_{11}}+\tfrac{1}{4}\big(\vec{v}_3\cdot\vec{1}-1\big)\right]\dd{U}{x}=\frac{v_{11}b_{0,1}-v_{12}b_{0,0}}{v_{11}-v_{12}}\,.\label{eq:paper1bc}
\end{equation}
This Robin boundary condition generates the accurate macroscale solution in the domain interior of Figure~\ref{fig:ch6num}.

To use the algorithm of Appendix~\ref{app:slowmanifoldbc}  to derive the left macroscale boundary conditions for system~\eqref{eq:sm} with different parameters, one needs to input the periodicity~$p$, the number of strands~$s$, the equilibrium spacing~$h$ and a $sp\times s(p+2)$ matrix summaries the first $sp$ equations in linear system~\eqref{eq:smss}.

\subsection{Right-end macroscale boundary condition\label{sub:First-iteration-1-1-1}}

Symmetry arguments give the right-end macroscale boundary condition.
Define a new longitudinal index~$n'=N-n$ but keep the same strand index~$j$.
The index~\(n'\) then parametrises the number of lattice points in from the right-hand boundary, just as the index~\(n\) does from the left-hand boundary.
We then repeat the  procedure described in Section~\ref{sub:left} with the new index~\(n'\).
To transform back to the original coordinate~$n$, the chain rule requires a change of the sign for the coefficient of~$\dd{U}{x}$, because of the change of direction of the longitudinal coordinate.

\subsection{Numerical results verify the two-strand two-periodic example\label{sec:ch7num2}}
This section numerically verifies the macroscale Robin boundary conditions discussed in Sections~\ref{sub:left} and~\ref{sub:First-iteration-1-1-1}, and confirms the validity of the quasi-steady state approximation.
Figure~\ref{fig:ch6num} plots the slowest eigenvector components of the microscale out-of-equilibrium spring-mass system~\eqref{eq:sm}.
We compare the slowest eigenvectors because they are most relevant to the large scale dynamics of the microscale system~\eqref{eq:sm}, and because the results were negligibly different for other cases.
As the system in Figure~\ref{fig:ch6num} only contains eight cells, the scale separations in this example is small.
Hence the assumption of infinite scale separation in homogenization theory \citep{MR2777986,PS08} is not appropriate here.
Nonetheless, our technique gives an accurate macroscale model.
For eigenvectors, the four microscale boundary values for this two-strand two-periodic system~$b_{0,0},b_{0,1},b_{N,0}$ and $b_{N,1}$ are zeros.
The blue line plots the slowest eigenvector of the  macroscale model~\eqref{eq:chp7m} with the classic macroscale Dirichlet boundary conditions~$U(0,t)=U(L,t)=0$\,.
The red line plots the slowest eigenvector of the macroscale model with derived Robin boundary conditions.
The triangles plot that of the microscale system~\eqref{eq:sm}.
Comparison between the two eigenvectors with different macroscale boundary values, for example at $10\leq n\leq 15$\,, demonstrates that the  derived boundary condition~\eqref{eq:paper1bc} improve the macroscale model~\eqref{eq:chp7m}.

The eigenvectors in Figure~\ref{fig:ch6num} reflect out-of-equilibrium dynamics.
These eigenvector plots show that the derived macroscale boundary condition~\eqref{eq:paper1bc} improves the macroscale model.
Hence the quasi-steady state assumption of Section~\ref{sub:left} is reasonable in this case.

\subsection{Numerical results verify a five-strand ten-periodic example\label{sec:ch7num}}

\begin{figure}
\centering\includegraphics{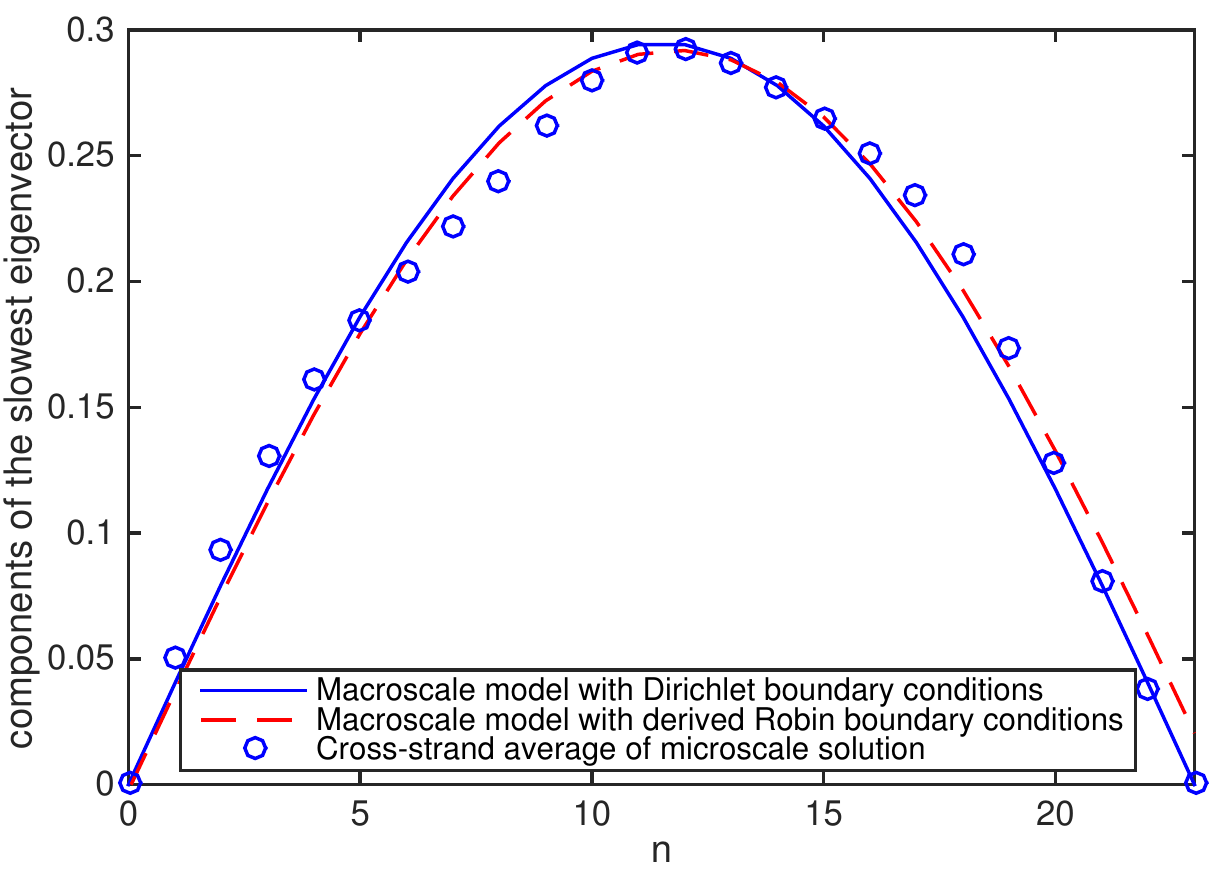}
\caption{The slowest eigenvector~$\vec{w}$ of a five-strand, ten-periodic, spring-mass system~\eqref{eq:sm} with microscale boundary values all zero.
The circles are the average of the slowest eigenvector of the microscale system with the same longitudinal index, that is,~$\tfrac{1}{5}\sum_{j=0}^{4}w_{n,j}$.
The red dotted line is our macroscale model~\eqref{eq:chp7m} obtained using boundary conditions~\eqref{eq:chp7bcf} that correctly produces the interior solution.
The blue line is the macroscale model obtained with classic heuristic zero boundary values, which is less accurate in the interior.\label{fig:ch7num}}
\end{figure}

We use the computer algebra algorithm to derive the macroscale boundary conditions for more complicated systems.
Figure~\ref{fig:ch7num} tests one five-strand ten-periodic example with~$N=23$\,.
Equations~\eqref{eq:ela} and~\eqref{eq:den} give the elasticity and density distributions of this example.
The ten microscale boundary values are all zero. To obtain the coefficient~$0.058$, one executes the function in Appendix~\ref{app:slowmanifoldmodel} by the command \texttt{bcfunction(A(1:s*p,1:s*(p+2)),s*p,h,p)}.
In this command, the matrix~$A$ is the $s(N-1)\times s(N+1)$ linear operator in the steady state system~\eqref{eq:smss}.
The macroscale boundary
conditions for this particular case are
\begin{equation}\label{eq:ch7ex}
U+0.058h\dd{U}{x}=0\quad\text{at }x=0\,,
\quad\text{and}\quad
U+0.53h\dd{U}{x}=0\quad\text{at }x=L.
\end{equation}
With these boundary conditions, the macroscale \pde~\eqref{eq:chp7m} has solutions that
accurately fit to the microscale solution within the interior of the
domain, better describing the global behaviour than the solution
with classic heuristic macroscale Dirichlet boundary conditions~\eqref{eq:paper1bc}.
As shown by the triangles in Figure~\ref{fig:ch7num}, the microscale boundary conditions at $x=23h$ force a boundary layer ($15\leq n\leq 23$) in the microscale model~\eqref{eq:sm}.
However, the macroscale \pde~\eqref{eq:chp7m} does not resolve the boundary layer.
Forcing the macroscale model to pass through the boundary layer  introduces an error in the interior of the domain, as shown by the solid blue line in Figure~\ref{fig:ch7num}.
Here we derive improved boundary conditions which reduces the interior error caused by poorly chosen macroscale boundary conditions.

\section{More complicated boundary conditions\label{sec:otherbc}}
Previous sections discuss microscale spring-mass
systems with Dirichlet boundary conditions.
This section considers different types of microscale boundary conditions
and determines the corresponding effective boundary conditions
for the spring-mass system~\eqref{eq:sm}.
Although the methodology in this section is general as in previous sections,  we use a two-strand two-periodic spring-mass system as an example for compactness.
Further, this section only treats the left-end boundary conditions: symmetry arguments, as in Section~\ref{sub:First-iteration-1-1-1}, correspondingly derive the right-end effective boundary conditions.

\subsection{Microscale specified flux boundary condition\label{sub:bcneumann}}

Specified flux Neumann boundary conditions are used in many mathematical models.
For example, petroleum engineers often model the dynamics of groundwater flow as diffusion in porous media with specified flux boundary conditions because the specified flux boundary condition best describes the physics of groundwater near the boundary surface \citep{WRCR:WRCR4207,WRCR:WRCR6440}.
Also, multiscale techniques are powerful in modelling porous media \citep{doi:10.1080/03605302.2014.926371,springerlink:10.1007/s00231-011-0879-4}.
Consequently, this section explores the effects of a microscale
specified flux boundary condition.

Suppose the spring-mass system~\eqref{eq:sm} is coupled
with specified flux boundary conditions at the left boundary of
\begin{equation}
u_{1,0}-u_{0,0}=hd_{0,0}\quad\text{and}\quad u_{1,1}-u_{0,1}=hd_{0,1}.\label{eq:bcneumannmicro}
\end{equation}
Analogous to equation~\eqref{eq:lincom}, write displacements~$\vec{u}_{0}$ as a linear combination of eigenvectors~$\vec{v}_{1}$,
$\vec{v}_{2}$, $\vec{v}_{3}$ and~$\vec{v}_{4}$ of a cell map~$T$
\begin{equation}
\vec{u}_{0}  =  \begin{bmatrix}u_{0,0}\\
u_{0,1}\\
u_{1,0}\\
u_{1,1}
\end{bmatrix}=c_{1}\vec{v}_{1}+c_{2}\vec{v}_{2}+c_{3}\vec{v}_{3}.\label{eq:otherbclin}
\end{equation}
Recall $\vec{v}_{2}=\left(1,1,1,1\right)$ is the eigenvector corresponding
to the eigenvalue of one.
Hence
\begin{equation}
\vec{u}_{0}=\begin{bmatrix}c_{1}v_{11}+c_{2}+c_{3}v_{31}\\
c_{1}v_{12}+c_{2}+c_{3}v_{32}\\
c_{1}v_{13}+c_{2}+c_{3}v_{33}\\
c_{1}v_{14}+c_{2}+c_{3}v_{34}
\end{bmatrix}.\label{eq:otherbcu0}
\end{equation}
In the case of Dirichlet boundary conditions, the microscale
boundary values~$u_{0,j},u_{N,j}$  and the definition of macroscale variables~\eqref{eq:amplitudeu} and~\eqref{eq:amplitudev}  force the
macroscale boundary values to satisfy constraint~\eqref{eq:twolayer-UdUdx}.
Analogous to the Dirichlet case, substitute equation~\eqref{eq:otherbcu0} into the microscale boundary constraint~\eqref{eq:bcneumannmicro} and relate
microscale and macroscale variables to obtain
\begin{equation}
\begin{bmatrix}v_{13}-v_{11} & 0 & v_{33}-v_{31}\\
v_{14}-v_{12} & 0 & v_{34}-v_{32}\\
0 & 1 & \tfrac{1}{4}\left(\vec{v}_{3}\cdot\vec{1}-1\right)\\
0 & 0 & \frac{1}{2h}
\end{bmatrix}\begin{bmatrix}c_{1}\\
c_{2}\\
c_{3}
\end{bmatrix}=\begin{bmatrix}hd_{0,0}\\
hd_{0,1}\\
\left.U\right|_{x=0}\\
\left.\dd{U}{x}\right|_{x=0}
\end{bmatrix}.\label{eq:twolayer-UdUdxneu2}
\end{equation}
The first two components of equation~\eqref{eq:twolayer-UdUdxneu2}
are constraints from the microscale boundary conditions~\eqref{eq:bcneumannmicro} and the last two
components are from the definition~\eqref{eq:amplitudeu} and~\eqref{eq:amplitudev}  of the macroscale
variable~$U$.

By inspection, the left null space of the \lhs~$4\times 3$ matrix in equation~\eqref{eq:twolayer-UdUdxneu2} is spanned by the vector
\begin{equation}\label{eq:nul1chp6}
\vec w = \begin{bmatrix} v_{14}-v_{12}
\\-v_{13}+v_{11}
\\0
\\2h\left[\left(v_{13}-v_{11}\right)\left(v_{34}-v_{32}\right)-\left(v_{14}-v_{12}\right)\left(v_{33}-v_{31}\right)\right] \end{bmatrix}.
\end{equation}
Premultiply constraint~\eqref{eq:twolayer-UdUdxneu2}
by vector~$\vec w ^T$ from~\eqref{eq:nul1chp6} to set the \lhs\ to zero.
Then constraint~\eqref{eq:twolayer-UdUdxneu2} gives
\begin{equation}
\left.\dd{U}{x}\right|_{x=0}=\frac{-\left(v_{14}-v_{12}\right)d_{0,0}+\left(v_{13}-v_{11}\right)d_{0,1}}{2\left[\left(v_{13}-v_{11}\right)\left(v_{34}-v_{32}\right)-\left(v_{14}-v_{12}\right)\left(v_{33}-v_{31}\right)\right]}.\label{eq:otherbcconstraint}
\end{equation}
The constraint~\eqref{eq:otherbcconstraint} forms a macroscale Neumann
boundary condition at $x=0$ for the macroscale \ode~\eqref{eq:chp7m}
corresponding to the microscale specified flux boundary condition~\eqref{eq:bcneumannmicro}.

\subsection{Microscale Robin-like boundary condition\label{sub:Microscale-'Robin-like'}}

In this section, suppose the spring-mass system~\eqref{eq:sm} has Robin-like boundary conditions at the left boundary of
\begin{equation}
u_{0,0}+\frac{d_{0,0}}{h}\left(u_{1,0}-u_{0,0}\right)=b_{0,0}\quad\text{and}\quad u_{0,1}+\frac{d_{0,1}}{h}\left(u_{1,1}-u_{0,1}\right)=b_{0,1}.\label{eq:bcrobinmicro}
\end{equation}
This boundary condition \eqref{eq:bcrobinmicro} is similar to, but
more complicated than, the specified flux boundary conditions~\eqref{eq:bcneumannmicro}.
Similar algebra as Section~\ref{sub:bcneumann} gives
\begin{equation}
\begin{bmatrix}v_{11}+\frac{1}{h}d_{0,0}(v_{13}-v_{11}) & v_{21}+\frac{1}{h}d_{0,0}(v_{23}-v_{21}) & v_{31}+\frac{1}{h}d_{0,0}(v_{33}-v_{31})\\
v_{12}+\frac{1}{h}d_{0,1}(v_{14}-v_{12}) & v_{22}+\frac{1}{h}d_{0,1}(v_{24}-v_{22}) & v_{32}+\frac{1}{h}d_{0,1}(v_{34}-v_{32})\\
0 & 1 & \tfrac{1}{4}\left(\vec{v}_{3}\cdot\vec{1}-1\right)\\
0 & 0 & \frac{1}{2h}
\end{bmatrix}\begin{bmatrix}c_{1}\\
c_{2}\\
c_{3}
\end{bmatrix}=\begin{bmatrix}b_{0,0}\\
b_{0,1}\\
\left.U\right|_{x=0}\\
\left.\dd{U}{x}\right|_{x=0}
\end{bmatrix}.\label{eq:twolayer-UdUdxrob}
\end{equation}
We need to compute a nullspace basis vector $\vec{w}=\left(w_{1},w_{2},w_{3},w_{4}\right)$ which satisfies
\[
\vec{w}^{T}\begin{bmatrix}v_{11}+\frac{1}{h}d_{0,0}(v_{13}-v_{11}) & v_{21}+\frac{1}{h}d_{0,0}(v_{23}-v_{21}) & v_{31}+\frac{1}{h}d_{0,0}(v_{33}-v_{31})\\
v_{12}+\frac{1}{h}d_{0,1}(v_{14}-v_{12}) & v_{22}+\frac{1}{h}d_{0,1}(v_{24}-v_{22}) & v_{32}+\frac{1}{h}d_{0,1}(v_{34}-v_{32})\\
0 & 1 & \tfrac{1}{4}\left(\vec{v}_{3}\cdot\vec{1}-1\right)\\
0 & 0 & \frac{1}{2h}
\end{bmatrix}=\vec{0}.
\]
Straightforward algebra gives
\begin{eqnarray*}
w_{1} & = & \frac{hv_{12}+d_{0,1}(v_{14}-v_{12})}{\left[hv_{11}+d_{0,0}(v_{13}-v_{11})\right]-\left[hv_{12}+d_{0,1}(v_{14}-v_{12})\right]},\\
w_{2} & = & -\frac{hv_{11}+d_{0,0}(v_{13}-v_{11})}{\left[hv_{11}+d_{0,0}(v_{13}-v_{11})\right]-\left[hv_{12}+d_{0,1}(v_{14}-v_{12})\right]},\\
w_{3} & = & 1,\\
w_{4} & = & -2h\frac{\left[hv_{12}+d_{0,1}(v_{14}-v_{12})\right]\left[hv_{31}+d_{0,0}\left(v_{33}-v_{31}\right)\right]}{\left[hv_{11}+d_{0,0}(v_{13}-v_{11})\right]-\left[hv_{12}+d_{0,1}(v_{14}-v_{12})\right]}\\
 &  & {}+2h\frac{\left[hv_{11}+d_{0,0}(v_{13}-v_{11})\right]\left[hv_{32}+d_{0,1}\left(v_{34}-v_{32}\right)\right]}{\left[hv_{11}+d_{0,0}(v_{13}-v_{11})\right]-\left[hv_{12}+d_{0,1}(v_{14}-v_{12})\right]}\\
 &  & {}-\tfrac{1}{2}h\left(\vec{v}_{3}\cdot\vec{1}-1\right).
\end{eqnarray*}
Premultiplying matrix-vector equation~\eqref{eq:twolayer-UdUdxrob} by vector~$\vec w^T$ gives the boundary condition at $x=0$ for the macroscale model~\eqref{eq:chp7m}
corresponding to the Robin-like microscale boundary condition~\eqref{eq:bcrobinmicro},
\[
U+w_{4}\dd{U}{x}=-w_{1}b_{0,0}-w_{2}b_{0,1}.
\]

\subsection{Microscale Cauchy-like boundary condition\label{sub:Microscale-'Cauchy-like'}}

For the microscale spring-mass system~\eqref{eq:sm}, suppose the specified boundary values are $u_{0,0}=b_{0,0}$ and $u_{1,0}=b_{1,0}$, that is,  fix the two leftmost values of strand~$j=0$ while leaving the two leftmost values~$u_{1,0}$ and~$u_{1,1}$ of strand~$j=1$ free to vary.
This is analogous to a Cauchy boundary condition\footnote{Recall that a Cauchy boundary condition in a continuous problem specifies the boundary value and the flux of the same boundary.} in a continuous problem because the flux in strand~$j=0$ is \((u_{1,0}-u_{0,0})/{h}\).
Knowing~$u_{0,0}$ and~$u_{1,0}$ is equivalent to knowing~$u_{0,0}$ and the flux.
Hence specifying function values~$u_{0,0}$ and~$u_{1,0}$ is similar to a Cauchy boundary condition, but only in one of the two strands.
Again, analogous to the matrix-vector equation~\eqref{eq:twolayer-UdUdx},
\begin{equation}
\begin{bmatrix}v_{11} & 1 & v_{31}\\
v_{13} & 1 & v_{33}\\
0 & 1 & \tfrac{1}{4}\left(\vec{v}_{3}\cdot\vec{1}-1\right)\\
0 & 0 & \frac{1}{2h}
\end{bmatrix}\begin{bmatrix}c_{1}\\
c_{2}\\
c_{3}
\end{bmatrix}=\begin{bmatrix}b_{0,0}\\
b_{1,0}\\
\left.U\right|_{x=0}\\
\left.\dd{U}{x}\right|_{x=0}
\end{bmatrix}.\label{eq:twolayer-UdUdxrob-1}
\end{equation}
Observe that this linear system is the same as system~\eqref{eq:twolayer-UdUdx} if we rename $v_{i3}$ to $v_{i2}$ for $i=1,2,3$ and $b_{1,0}$ to $b_{0,1}$.
Hence applying  mapping~$v_{i2}\mapsto v_{i3}$ for $i=1,2,3$ and $u_{0,1} \mapsto b_{1,0}$ to the macroscale boundary conditions~\eqref{eq:paper1bc} gives a Robin macroscale boundary condition for the macroscale model~\eqref{eq:chp7m} at $x=0$ corresponding to Cauchy-like microscale boundary condition,
\[
U-2h\left[\frac{v_{11}v_{33}-v_{13}v_{31}}{v_{13}-v_{11}}+\tfrac{1}{4}\left(\vec{v}_3\cdot\vec{1}-1\right)\right]\dd{U}{x}=\frac{v_{11}b_{1,0}-v_{13}b_{0,0}}{v_{11}-v_{13}}\,.
\]

\subsection{Microscale mixed boundary conditions}

The procedures used to model Robin-like microscale boundary conditions
\eqref{eq:bcrobinmicro} and Cauchy-like microscale boundary conditions in Section~\ref{sub:Microscale-'Cauchy-like'}
are almost the same as the procedure in Section~\ref{sub:bcneumann} for specified flux microscale boundary conditions.
However, a combination of these microscale boundary conditions described in Section~\ref{sub:bcneumann},
\ref{sub:Microscale-'Robin-like'} and \ref{sub:Microscale-'Cauchy-like'} leads
to a very different situation.

Suppose the two-strand spring-mass system~\eqref{eq:sm}
has mixed boundary conditions.
On the left boundary, we specify
three boundary values as constants~$u_{0,0}=b_{0,0},\ u_{0,1}=b_{0,1}$ and $u_{1,0}=b_{1,0}$\,.
This microscale boundary condition is equivalent to specifying two boundary values and the flux between~$u_{0,0}$ and~$u_{1,0}$.
On the right boundary, we specify just one condition, the boundary value of the last field on the top strand~$u_{N,0}=b_{N,0}$\,.
We assume the cross elasticity near the right boundary
is not small enough for~$u_{N,1}$ to be unbounded.
This assumption
is important because if~$u_{N,1}$ was unbounded, then the unstable mode of cell map~$T$ in the spatial evolution would be non-trivial.

On the left boundary, we have one more specified boundary value~$u_{0,1}$
than the Cauchy-like boundary condition analysed in Section~\ref{sub:Microscale-'Cauchy-like'}.
Hence we adjoin the boundary value~$u_{0,1}$ on the bottom strand,
that is, the second component of matrix-vector equation~\eqref{eq:twolayer-UdUdx},
to matrix-vector equation~\eqref{eq:twolayer-UdUdxrob-1},

\begin{equation}
\begin{bmatrix}v_{11} & 1 & v_{31}\\
v_{13} & 1 & v_{33}\\
v_{12} & 1 & v_{32}\\
0 & 1 & \tfrac{1}{4}\left(\vec{v}_{3}\cdot\vec{1}-1\right)\\
0 & 0 & \frac{1}{2h}
\end{bmatrix}\begin{bmatrix}c_{1}\\
c_{2}\\
c_{3}
\end{bmatrix}=\begin{bmatrix}b_{0,0}\\
b_{1,0}\\
b_{0,1}\\
\left.U\right|_{x=0}\\
\left.\dd{U}{x}\right|_{x=0}
\end{bmatrix}.\label{eq:twolayer-UdUdxrob-2}
\end{equation}
Recall that this matrix-vector equation is solvable only if the \rhs\ is in the range, that is, orthogonal to the left null space of the \lhs\ $5\times3$ matrix.
By direct substitution one verifies the two linear independent vectors,
\[
\vec{w}_{1}=\begin{bmatrix}
v_{13}/(v_{11}-v_{13})\\
-v_{11}/\left(v_{11}-v_{13}\right)\\
0\\
1\\
-2h\left(\frac{v_{13}v_{31}}{v_{11}-v_{13}}-\frac{v_{11}v_{33}}{v_{11}-v_{13}}+\tfrac{1}{4}\left(\vec{v}_{3}\cdot\vec{1}-1\right)\right)
\end{bmatrix}
\]
and
\[
\vec{w}_{2}=\begin{bmatrix}v_{12}/(v_{11}-v_{12})\\
0\\
-v_{11}/(v_{11}-v_{12})\\
1\\
-2h\left(\frac{v_{12}v_{31}}{v_{11}-v_{12}}-\frac{v_{11}v_{32}}{v_{11}-v_{12}}+\tfrac{1}{4}\left(\vec{v}_{3}\cdot\vec{1}-1\right)\right)
\end{bmatrix},
\]
span the two dimensional null space of the matrix on the \lhs\
of boundary constraint~\eqref{eq:twolayer-UdUdxrob-2}.
Premultiplying
the constraints~\eqref{eq:twolayer-UdUdxrob-2} by either~$\vec{w}_{1}^{T}$
or~$\vec{w}_{2}^{T}$ zeros the \lhs.
Then rearrange the expression on the \rhs,
\begin{eqnarray*}
U-2h\left[\frac{v_{13}v_{31}-v_{11}v_{33}}{v_{11}-v_{13}}+\tfrac{1}{4}\left(\vec{v}_3\cdot\vec{1}-1\right)\right]\dd{U}{x}&=&\frac{v_{11}b_{1,0}-v_{13}b_{0,0}}{v_{11}-v_{13}}\, ,\nonumber \\
U-2h\left[\frac{v_{12}v_{31}-v_{11}v_{32}}{v_{11}-v_{12}}+\tfrac{1}{4}\left(\vec{v}_3\cdot\vec{1}-1\right)\right]\dd{U}{x} &= &\frac{v_{11}b_{0,1}-v_{12}b_{0,0}}{v_{11}-v_{12}}\,.
\end{eqnarray*}
Combining these two boundary constraints in two different ways gives
a Cauchy boundary condition for the macroscale model~\eqref{eq:chp7m}
at $x=0$\,:
\begin{eqnarray}\label{eq:mbcos}
U & = & \frac{v_{11}b_{1,0}-v_{13}b_{0,0}}{v_{11}-v_{13}}  -  \left[\frac{v_{13}v_{31}-v_{11}v_{33}}{v_{11}-v_{13}}+\tfrac{1}{4}\left(\vec{v}_3\cdot\vec{1}-1\right)\right]\mathcal{D},\nonumber\\
\dd{U}{x} & = & -\frac{\mathcal{D}}{2h},
\end{eqnarray}
where
\[
\mathcal{D}=\frac{v_{11}\left[\left(v_{13}-v_{11}\right)\left(u_{0,1}-u_{0,0}\right)-\left(v_{12}-v_{11}\right)\left(b_{1,0}-b_{0,0}\right)\right]}{\left(v_{11}v_{33}-v_{31}v_{13}\right)\left(v_{12}-v_{11}\right)+\left(v_{12}v_{31}-v_{11}v_{32}\right)\left(v_{13}-v_{11}\right)}.
\]
Recall that the macroscale model~\eqref{eq:chp7m} is a second order \pde , hence two boundary conditions at the left-hand end is enough for the macroscale model.

No more boundary conditions are needed at the right-hand end because the second order macroscale wave \pde ~\eqref{eq:chp7m} only need two boundary conditions. Interestingly, the microscale boundary value~$b_{N,0}$ is not present in the macroscale boundary conditions~\eqref{eq:mbcos} and gives no macroscale constraint.

\section{Conclusion\label{sec:ch7sum}}
This article generalises and extends methodologies piloted by \citet{Chen2014} to $s$-strand $p$-periodic spring-mass systems.
Section~\ref{sec:ch7int} finds the macroscale homogeneous model for the microscale inhomogeneous microscale spring-mass system~\eqref{eq:sm} with Dirichlet boundary conditions using computer algebra.
Section~\ref{sub:First-iteration-1ch7} computes the macroscale Robin boundary conditions for the macroscale model.
Section~\ref{sec:otherbc} generalises the methodology of deriving effective macroscale boundary conditions to spring-mass systems with more complicated microscale boundary conditions.
Section~\ref{sec:ch7num} verifies that the derived boundary conditions improves the accuracy of the macroscale model.

As we did not use the hyperbolicity of the spring-mass system, our methodology is applicable to analogous heterogeneous diffusion systems in long-thin domain.

The proposed approach to find boundary conditions for macroscale models is currently only applicable to problems in long thin domains.
In a general multi-dimensional problem, the stable manifold near the boundary usually decays differently in different directions.
The derivation of macroscale boundary conditions in the discussed problems involves finding an operator that maps the solution in one cell to the next cell.
In cases where a domain which is large in two or three dimensions, it is hard to find an equivalent of ``the next cell''.
Extending this methodology to multi-dimensions requires us to find a scheme that parametrises the decaying stable manifold near boundaries with multiple macroscale directions.

Moreover, this article assumes the time derivative terms have negligible impact on the boundary structure when deriving macroscale boundary conditions.
The effect of time dependence upon macroscale boundary conditions is a subject for further research.

\bibliographystyle{agsm}
\bibliography{research}

\input{appendix}

\end{document}

%% file: appendix.tex
\appendix
\section{Coefficients for slow manifold\label{app:slowmanifold}}
This appendix lists the coefficients~$\vec {\alpha}$ and $\vec {\beta}$ of the slow manifold~\eqref{eq:ch7ex}:
\begin{eqnarray*}
\vec {\alpha}&=&(
-0.5186	,
0.001565	,
-0.05784	,
0.1854	,
0.2666	,
0.0335	,
0.4125	,
0.3109	,
0.32	,
0.4105	\\&&
0.3952	,
0.7044	,
0.5262	,
0.2587	,
0.3662	,
0.6026	,
0.8001	,
0.5527	,
0.1362	,
0.2647	\\&&
0.6325	,
0.6093	,
0.3736	,
0.03515	,
0.1916	,
0.6048	,
0.09476	,
0.0295	,
-0.08023	,
0.02459	\\&&
0.3704	,
-0.4336	,
-0.3396	,
-0.24		,
-0.2676	,
-0.3766	,
-0.7119	,
-0.5954	,
-0.3766	,
-0.3865	\\&&
-1.055	,
-0.6954	,
-0.652	,
-0.3558	,
-0.3074	,
-1.046	,
-0.4202	,
-0.4465	,
-0.1149	,
-0.03636	);\\
\vec {\beta}&=&(-0.59	,
0.3042	,
-0.6824	,
-0.8	,
-0.7215	,
0.6218	,
0.5918	,
-0.04535	,
-0.3774	,
-0.6416	\\&&
1.1	,
0.4504	,
0.4774	,
0.09984	,
-0.3785	,
0.7149	,
0.1631	,
0.7654	,
0.4951	,
-0.09197	\\&&
0.09442	,
0.1039	,
0.906	,
0.7926	,
0.2169	,
-0.5301	,
0.2259	,
0.8594	,
0.9142	,
0.4569	\\&&
-0.7361	,
0.09054	,
0.5265	,
0.7474	,
0.6825	,
-0.8021	,
-0.2024	,
-0.05982	,
0.2586	,
0.6558	\\&&
-1.392	,
-0.348	,
-0.6996	,
-0.4206	,
0.1816	,
-1.539	,
-0.1533	,
-0.9945	,
-0.8702	,
-0.4211)	
\end{eqnarray*}

\section{Computer algebra derives slow manifold model \label{app:slowmanifoldmodel}}
This appendix lists the computer algebra code which derives the slow manifold model in Section~\ref{sub:slowm}.
The code is written in Reduce which is a free general purpose computer algebra package.\footnote{\protect\url{http://www.reduce-algebra.com/}}
\citet{Fateman:2003:CSP:844076.844080} demonstrated that Reduce is more than twenty times faster than the Matlab symbolic toolbox (which invokes Mupad) and generally over an order of magnitude faster than other popular packages.


\input{code/wstrand_Nstep}

\section{Computer algebra derives the macroscale boundary conditions \label{app:slowmanifoldbc}}


This appendix lists the Matlab code which derives the macroscale boundary conditions in Section~\ref{sub:First-iteration-1ch7}.


\begin{lstlisting}
function [coe_dudx] = bcfunction(diff,sp,h,p)
%{
Compute boundary condition U+d*dudx=b chen chen 01/09/2014
This algorithm derives the macroscale boundary condition. It
computes a vector coe_dudx of size $s+2$. The coefficients
of derived Robin boundary conditions are coponents of
coe_dudx This algorithm also finds the macroscale boundary
conditions at the right boundary x=L using symmetry
arguments.

The input matrix $diff$ consists of the first sp rows of the
quasi steady state spring-mass system in matrix-vector form.
As the dynamics of the s left-most masses and the s
right-most on each masses within a cell are coupled with s
masses closest to the cell side, the matrix diff is sp
\times s(p+2). If matrix diff is of data type double, Matlab
derives boundary coefficients~$d_{0}$ and $d_{L}$
numerically. If matrix diff is a symbolic matrix, Matlab
invokes the Symbolic Toolbox to derive symbolically the
boundary condition coefficients.

The input sp is the total number of masses within a cell.
The input p is the periodicity of the spring-mass system.
The input h is the separation between two adjacent lattice
points. s is the number of strand
%}
s=sp/p;

%{
Compute the $2s\times 2s$ square matrix map which maps the
displacement of the first 2s masses in a cell to that of the
left-most 2s masses in the next cell. Also compute the
$(p-2)s\times2s$ matrix map_cell1 which maps the first 2s
masses in a cell to the remaining (p-2)s masses within the
same cell. Two-periodic systems have the matrix map_cell1
empty because the longitudinal diffusivities alternate in
the two-periodic systems. Matrix map_cell1 is only used to
compute the macroscale varaibles U and its derivative dUdx.
%}
map_whole=-diff(:,s*2+1:sp+2*s)\diff(:,1:s*2);
map_cell1=map_whole(1:sp-2*s,:);
map=map_whole(end-2*s+1:end,:);

[evecn evaln]=eig(map);
%discard small numerical errors
j=find(abs(imag(evecn))<1e-6); evecn(j)=real(evecn(j));
j=find(abs(imag(evaln))<1e-6); evaln(j)=real(evaln(j));
diagevaln=diag(evaln); % diagonalise to be ordered
if sum(real(diagevaln)<=0)~=0
    disp('negative eigenvalue');
    pause;
end
[~,order]=sort(-diagevaln,'descend');%
evaln=(diag(diagevaln(order))) %reorder
evecn=(evecn(:,order)); % reorder

%compute generalised eigenvector
v1n=evecn(:,sp/p);
%1 is the corresponding eigenvalue
gv=([map-1*eye(2*s) -v1n;-v1n' 0])\rand(2*s+1,1);
gv=gv(1:2*s)/gv(end);
% put the generalised eigenvector into eigenvector matrix
evecn(:,sp/p+1)=gv; 

u0=evecn 
u1=map_cell1*u0 
U_cell1=mean([u0;u1]) %rows are U(1/2),U(1/2+p)
if sum(v1n/mean(v1n)-ones(2*s,1))<1e-7
    dUdx=mean(v1n)/p/h;
else
    disp('eigenvector not ones');
    pause;
end
%[unstable centre centre stable]
dUdx=[zeros(1,s-1) 0 dUdx zeros(1,s-1)]; 
U0=U_cell1-dUdx*mean(0:p-1)*h
UdUdx=[U0;dUdx]

bcuux=[evecn(1:s,:);UdUdx]
bcuux(s+1:end,1:(s-1))=zeros(2,s-1);

%{
Find the vector zz, which contains the coefficients of the
derived boundary conditions, by finding the null space of a
$(s+2)\times(s+1)$ linear system. 
%}
zz=null(bcuux(:,1:(s+1))');

%{
normalise by put the coefficient of U as 1. The macroscale
boundary condition coefficient $d_{0}$ in the derived Robin
macroscale boundary conditions is the last component of
vector coe_dudx. The Robin boundary value~$B_{0}(t)$ in the
derived macroscale boundary conditions is the dot product
$B_{0}(t)=-(u_{0,0},u_{0,1},\ldots,u_{0,s-1},0,0)\cdot
coe_dudx$ where~$u_{0,0},u_{0,1},\ldots,u_{0,s-1}$ are the
specified Dirichlet boundary values for the quasi steady
state wave equation.
%}
coe_dudx=zz(end)/zz(end-1)/h %coefficient of h*dUdx

%clear small numerical error
j=find(abs(imag(coe_dudx))<1e-6); coe_dudx(j)=real(coe_dudx(j));
end
\end{lstlisting}


%% file: code/wstrand_Nstep.tex
Use iteration to form the slow manifold model of the s-strand
 p-periodic wave pde. Here the description is in terms of
 the mean displacement \(U\) and the mean velocity \(V\) over a cell.
 \(U\) and \(V\) are in Fourier space. (Chen Chen, August 2014)

Formatting for printed output
\begin{reduce}
on div; off allfac; on revpri; factor k,u,v;
on rounded; print_precision 4$
\end{reduce}
provide linear algebra functions
\begin{reduce}
load_package linalg;
\end{reduce}
 Firstly, define the macroscale parametric variables $\tilde{U}$ and
$\tilde{V}$ to depend upon time; and secondly define time derivatives of $\tilde{U}$ and $\tilde{V}$, that is $\dd{\tilde{U}}{t}$
and~$\dd{\tilde{V}}{t}$, as \texttt{V}
and~\texttt{g}, since $\tilde V=\dd{\tilde{U}}{t}$ by the amplitude definitions.
and~\texttt{g} stores the current approximation of the macroscale evolution.

\begin{reduce}
depend U,t;let df(U,t)=>V; let df(V,t)=>g;
\end{reduce}

Define a function which computes the modulus of two integers. 
This routine constructs matrix~$\mathcal{L}_{k}$ to errors~$\Ord(k^3)$.
\begin{reduce}
procedure mod(ii,nos);
    round((ii/nos-floor(ii/nos))*nos);
\end{reduce}

\paragraph{Input parameters}
Read the parameters from a file named \texttt{parameter.txt}.
Within this file, only the last three lines of codes
are compulsory. File \texttt{parameter.txt} has to specify the
number of strands and the periodicity of the spring-mass system.
\begin{reduce}
in "parameters.txt";
\end{reduce}
The number of masses in a cell
\begin{reduce}
non:=nos*nol;
\end{reduce}

\paragraph{Set initial approximation}
Assign the initial linear slow manifold approximation~U
and evolution approximation~g to variables~\texttt{vv} and~\texttt{uu}. These
lines also stores the matrix~$B$ in fourier spaced microscale system
 in the variable~\texttt{rhom}.

\begin{reduce}
matrix rhom(non,non),vv(non,1),uu(non,1),ones(non,1),zero(non,1);
g:=0;count:=0;
for ii:=1:nos do
  for j:=1:nol do begin
    count:=count+1;
    rhom(count,count):=mkid(mkid(rho,ii-1),j-1)*h^2;
    vv(count,1):=V;
    uu(count,1):=U;
    ones(count,1):=1;
    zero(count,1):=1;
  end;
clear count;
\end{reduce}

\paragraph{Code the wave dynamical system}
 Forms matrix \(L_k\). Only keep terms up to order two in \(k\). Edit this line for higher order truncation
Taylor expand $e^{\pm ikh}$ as $1\pm ikh-\frac{1}{2}h^2k^{2}+\mathcal{O}(k^{3})$.
The \texttt{let} command controls the truncation of the asymptotic
approximation by informing Reduce to always discard any factor in $k^{3}$
or higher. Thus all expressions are computed to an error of $\Ord(k^{3})$.
One advantage of Reduce algebra is this `\texttt{let}' command. 
 With a single command \verb|let k^3=>0|,
I discard any factor in $k^{3}$ or higher.

\begin{reduce}
let {k^3=>0};
expkpos:=1+i*h*k-1/2*h^2*k^2;expkneg:=1-i*h*k-1/2*h^2*k^2;
\end{reduce}

Set up matrix~$\mathcal{L}_{k}$ to errors $\mathcal{O}(k^{3})$.
\begin{reduce}
array kk(nos,nol,nol),kappa(nos,nol,nol);
for ii:=0:nos-1 do
  for j:=0:nol-1 do begin
    kk(ii,j+1,j+1):=mkid(mkid(k,ii),j);
    kappa(ii,j+1,j+1):=kappa(ii,j+1,j+1)
    -mkid(mkid(k,ii),j)-mkid(mkid(k,mod(ii-1,nos)),j);    
    for jj:=0:nol-1 do
      if neq(jj,j) then <<
      kappa(ii,j+1,jj+1):=
        mkid(mkid(mkid(kc,ii),min(jj,j)),max(jj,j));    
      kappa(ii,j+1,j+1):=kappa(ii,j+1,j+1)
        -mkid(mkid(mkid(kc,ii),min(jj,j)),max(jj,j));
        >>;    
  end;
matrix Lk(non,non);
for ii:=0:nos-1 do
  for j1:=0:nol-1 do begin
    lk(j1+1+ii*nol,mod(j1+ii*nol+nol,non)+1)
    :=lk(j1+1+ii*nol,mod(j1+ii*nol+nol,non)+1)
    +kk(ii,j1+1,j1+1)*expkpos;
    lk(j1+1+ii*nol,mod(j1+ii*nol-nol,non)+1)
    :=lk(j1+1+ii*nol,mod(j1+ii*nol-nol,non)+1)
    +kk(mod(ii-1,nos),j1+1,j1+1)*expkneg;
    for j2:=0:nol-1 do
      lk(j1+1+ii*nol,j2+1+ii*nol):=kappa(ii,j1+1,j2+1);
    end;
  end;
\end{reduce}

Compute matrix~$\mathcal{L}_{0}$ by substituting $k=0$ into matrix~$\mathcal{L}_{k}$.
\begin{reduce}
L0:=sub(k=0,Lk);
\end{reduce}

Constrain the sum of corrections~\verb|uuhat| and \verb|vvhat|
to zero by replacing two of the update equations. 
By doing this, the amplitude definition are satisfied. 
We can replace two of the equations because the matrix~$\mathcal L_0$ has 
eigenvalues of zero with multiplicity of two. 
By adjoining two amplitude equations 
the corrections are still  unique.
\begin{reduce}
temp:=L0;
for j:=1:nol*nos do temp(nol*nos,j):=1;
\end{reduce}

\paragraph{Iteratively construct centre manifold model}
\begin{reduce}
for iter:=1:12 do begin
\end{reduce}
    Compute the residual of the current approximation.
\begin{reduce}
  write res1:=df(uu,t)-vv;
  write res2:=rhom*df(vv,t)-Lk*uu;
\end{reduce}
    This step is the solvability condition.
    The solvability condition chooses the correction of 
    the evolution $\dd{V}{t}$ so that the slow manifold 
    corrections uuhat and vvhat exist.
\begin{reduce}
  ghat:=TP(ones)*res2/(-for ii:=1:nos*nol sum rhom(ii,ii));
  write g:=g+ghat(1,1);
\end{reduce}
    Solve for the correction of the slow manifold model uuhat
    and vvhat and update current approximation.
\begin{reduce}
  tempv:=res2+rhom*ghat(1,1)*ones;
  tempv(nos*nol,1):=0;
  uuhat:=TP(TP(tempv)/TP(temp));
  uuhat(nos*nol,1):=-for ii:=1:nos*nol-1 sum uuhat(ii,1);
  uu:=uu+uuhat;
  vvhat:=res1+df(uuhat,U)*V;
  vvhat(nos*nol,1):=-for ii:=1:nos*nol-1 sum vvhat(ii,1);
  vv:=vv+vvhat;
  showtime;
\end{reduce}
    Terminate the iteration when both residuals res1 and res2 are \(\Ord(k^3)\).
\begin{reduce}
  if {res2,res1}={zero,zero} then write iter:=1000000+iter;
end;
\end{reduce}
Output the macroscale effective coefficients~$\bar{\kappa}/\bar{\rho}$
to a file. Printing to file provides a means for Reduce computer algebra
to communicate with any numerical scheme.
write output for numerics to read
\begin{reduce}
out macro;
write -coeffn(coeffn(g,k,2),u,1)$
shut macro$
end;
\end{reduce}